\documentclass[12pt]{article}

\usepackage{indentfirst}
\usepackage{amsfonts}
\usepackage{amsmath}
\usepackage{amssymb}
\usepackage{amsbsy, amsthm}
\usepackage[margin=1in]{geometry}

\newtheorem{dfn}{Definition}[section]
\newtheorem{theorem}[dfn]{Theorem}
\newtheorem{lemma}[dfn]{Lemma}
\newtheorem{prop}[dfn]{Proposition}

\newtheorem{corollary}[dfn]{Corollary}
\newtheorem{question}[dfn]{Question}
\newenvironment{pf}{\noindent{\bf Proof.}}
{\enspace\vrule height5pt depth0pt width5pt} 

\def\F {{\mathcal F}}

\def\P {{\mathcal P}}
\def\X {{\mathcal X}}
\def\W {{\mathcal W}}

\def\U {{\mathcal U}}

\def\ad {{\rm asdim}}

\begin{document}

\title{Asymptotic dimension of minor-closed families and beyond}
\author{Chun-Hung Liu\thanks{chliu@math.tamu.edu. Partially supported by NSF under Grant No.~DMS-1929851 and DMS-1954054.} \\
\small Department of Mathematics, \\
\small Texas A\&M University,\\
\small College Station, TX 77843-3368, USA \\
}

\maketitle

\begin{abstract}
The asymptotic dimension of metric spaces is an important notion in geometric group theory introduced by Gromov.
The metric spaces considered in this paper are the ones whose underlying spaces are the vertex-sets of graphs and whose metrics are the distance functions in graphs.
A standard compactness argument shows that it suffices to consider the asymptotic dimension of classes of finite graphs.

In this paper we prove that the asymptotic dimension of any proper minor-closed family, any class of graphs of bounded tree-width, and any class of graphs of bounded layered tree-width are at most 2, 1, and 2, respectively.
The first result solves a question of Fujiwara and Papasoglu; the second and third results solve a number of questions of Bonamy, Bousquet, Esperet, Groenland, Pirot and Scott.
These bounds for asymptotic dimension are optimal and improve a number of results in the literature.
Our proofs can be transformed into linear or quadratic time algorithms for finding coverings witnessing the asymptotic dimension which is equivalent to finding weak diameter colorings for graphs.
The key ingredient of our proof is a unified machinery about the asymptotic dimension of classes of graphs that have tree-decompositions of bounded adhesion over hereditary classes with known asymptotic dimension, which might be of independent interest.
\end{abstract}

\section{Introduction} \label{sec:intro}

The asymptotic dimension of metric spaces is an important notion in geometric group theory introduced by Gromov \cite{g}.
(See \cite{bd} for a survey.)
There are a number of equivalent definitions for asymptotic dimension.
Here we use the following version.
For a nonnegative integer $n$, the {\it $n$-dimensional control function of a metric space} $(X,d)$ is a function $f: {\mathbb R}^+ \rightarrow {\mathbb R}^+$ such that for every positive real number $r$, there exist collections $\U_1,\U_2,...,\U_{n+1}$ of subsets of $X$ such that 
	\begin{itemize}
		\item $\bigcup_{i=1}^{n+1}\bigcup_{U \in \U_i}U \supseteq X$,
		\item for each $i \in [n+1]$, if $U,U'$ are distinct elements of $\U_i$, then $d(x,x')>r$ for every $x \in U$ and $x' \in U'$, and
		\item for each $i \in [n+1]$ and $U \in \U_i$, if $x,x' \in U$, then $d(x,x') \leq f(r)$.
	\end{itemize}
The {\it asymptotic dimension of a metric space} $(X,d)$ is the least integer $n$ such that there exists an $n$-dimensional control function.

For example, when $(X,d)=({\mathbb R}, \| \cdot \|_2)$, $f(x):=x$ is a 1-dimensional control function of $(X,d)$, since for every positive real number $r$, we may take $\U_1 = \{[mr,(m+1)r): m$ is an odd integer$\}$ and $\U_2 = \{[mr,(m+1)r): m$ is an even integer$\}$ so that the above conditions are satisfied.
Similarly, one can show that $({\mathbb R}^2, \| \cdot \|_2)$ has asymptotic dimension at most 2 by first dividing ${\mathbb R}^2$ into regular hexagons with side length $2r$ and then putting then into $\U_1,\U_2,\U_3$ alternately.

In this paper we study the asymptotic dimension of the metric spaces defined by graphs.
For a nonnegative integer $n$, the {\it $n$-dimensional control function of a graph} $G$ is an $n$-dimensional control function of the metric space $(V(G),d_G)$, where for every $x,y \in V(G)$, $d_G(x,y)$ is the minimum number of edges of a path in $G$ between $x$ and $y$.
We call $d_G(x,y)$ the {\it distance in $G$ between $x$ and $y$}.
The {\it asymptotic dimension of a graph} $G$ is the minimum $n$ such that there exists an $n$-dimensional control function of $G$. 
For a class $\F$ of graphs, the {\it asymptotic dimension} of $\F$, denoted by $\ad(\F)$, is the minimum $n$ such that there exists a common $n$-dimensional control function of all graphs in $\F$.

Note that the asymptotic dimension of any finite graph is 0.
So it is only interesting to consider asymptotic dimension of an infinite graph or an infinite class of graphs.
A simple compactness argument shows that the asymptotic dimension of an infinite graph is at most the asymptotic dimension of the class of its finite induced subgraphs.

Hence in this paper, we only consider the asymptotic dimension of classes of finite graphs.
From now on, {\it graphs} in this paper are finite, unless otherwise stated.

The work of this paper is motivated by the following question of Fujiwara and Papasoglu \cite{fp} about asymptotic dimension and graph minors.
A graph $H$ is a {\it minor} of another graph $G$ if $H$ is isomorphic to a graph that can be obtained from a subgraph of $G$ by contracting edges.

\begin{question}[{\cite[Question 5.2]{fp}}] \label{question_minor}
Is there a constant $M$ such that for every graph $H$, the class of $H$-minor free graphs has asymptotic dimension at most $M$?
Can we take $M=2$?
\end{question}

Ostrovskii and Rosenthal \cite{or} proved a bound depending on the graph $H$: for every graph $H$, the class of $H$-minor free graphs has asymptotic dimension at most $4^{\lvert V(H) \rvert}-1$.
Bonamy, Bousquet, Esperet, Groenland, Pirot and Scott \cite{bbegps} proved the case for bounded maximum degree graphs: for every integer $\Delta$ and graph $H$, the class of $H$-minor free graphs of maximum degree at most $\Delta$ has asymptotic dimension at most 2.
Bonamy et al.\ \cite{bbegps} used their result to answer a question of Ostrovskii and Rosenthal \cite{or} by showing that for every finitely generated group $\Gamma$ with a symmetric finite generating set $S$, if there exists a graph $H$ such that $H$ is not a minor of the Cayley graph for $(\Gamma,S)$, then the asymptotic dimension of the Cayley graph for $(\Gamma,S)$ is at most 2.
Bonamy et al.\ \cite{bbegps} also proved the case $H=K_{3,t}$ for every positive integer $t$ and the case that $H$ is an apex-forest in Question \ref{question_minor}.

One of the main results of this paper is a complete solution of Question \ref{question_minor} and hence improves the aforementioned results of Ostrovskii and Rosenthal \cite{or} and Bonamy et.\ al.\ \cite{bbegps}.
It also solves a question \cite[Question 5]{bbegps} about forbidding $K_{s,t}$-minors in a strong sense.

\begin{theorem} \label{minor_ad_intro}
For every graph $H$, the asymptotic dimension of the class of $H$-minor free graphs is at most 2.
\end{theorem}

Note that the number 2 in Theorem \ref{minor_ad_intro} is optimal if $H$ is a non-planar graph, since the class of 2-dimensional grids has asymptotic dimension 2 \cite{g}.
When $H$ is planar, we prove that the number 2 can be improved to be 1.

\begin{theorem} \label{tw_ad_intro_1}
For every planar graph $H$, the asymptotic dimension of the class of $H$-minor free graphs is at most 1.
\end{theorem}

Bonamy et al.\ \cite[Question 3]{bbegps} asked whether there exists a constant $k$ such that for every positive integer $w$, the class of graphs of tree-width at most $w$ has asymptotic dimension at most $k$.
A {\it tree-decomposition} of a graph $G$ is a pair $(T,\X)$ such that $T$ is a tree and $\X$ is a collection $(X_t: t \in V(T))$ of subsets of $V(G)$, called the {\it bags}, such that
	\begin{itemize}
		\item $\bigcup_{t \in V(T)}X_t = V(G)$,
		\item for every $e \in E(G)$, there exists $t \in V(T)$ such that $X_t$ contains the ends of $e$, and
		\item for every $v \in V(G)$, the set $\{t \in V(T): v \in X_t\}$ induces a connected subgraph of $T$.
	\end{itemize}
For a tree-decomposition $(T,\X)$, the {\it adhesion} of $(T,\X)$ is $\max_{tt' \in E(T)}\lvert X_t \cap X_{t'} \rvert$, and the {\it width} of $(T,\X)$ is $\max_{t \in V(T)}\lvert X_t \rvert-1$.
The {\it tree-width} of $G$ is the minimum width of a tree-decomposition of $G$.

Note that by the Grid Minor Theorem \cite{rs_V}, excluding any planar graph as a minor is equivalent to having bounded tree-width.
Hence Theorem \ref{tw_ad_intro_2} which is an equivalent form of Theorem \ref{tw_ad_intro_1} solves this question of Bonamy et al.\ in a strong sense.

\begin{theorem} \label{tw_ad_intro_2}
For every positive integer $w$, the asymptotic dimension of the class of graphs of tree-width at most $w$ is at most 1.
\end{theorem}

Bell and Dranishnikov \cite{bd} and Fujiwara and Papasoglu \cite{fp}, respectively, showed that the class of trees and the class of cacti, respectively, have asymptotic dimension 1.
Note that trees have tree-width at most 1 and cacti have tree-width at most 2.
So Theorems \ref{tw_ad_intro_1} and \ref{tw_ad_intro_2} are optimal and generalize the results of Bell and Dranishnikov \cite{bd} and Fujiwara and Papasoglu \cite{fp}.
Theorems \ref{tw_ad_intro_1} and \ref{tw_ad_intro_2} also generalize results of Bonamy et al.\ \cite{bbegps} who proved the same for the class of graphs of bounded path-width and proved the case with the extra bounded maximum degree condition, where the latter was implicitly proved in an earlier paper of Benjamini, Schramm and Tim\'{a}r \cite{bst}.

By using a trick in \cite{bbegps}, we are able to obtain optimal results for asymptotic dimension for graphs of bounded layered tree-width.
A {\it layering} of a graph $G$ is an ordered partition $(V_1,V_2,...)$ of $V(G)$ into (possibly empty) subsets such that for every edge $e$ of $G$, there exists $i_e$ such that $V_{i_e} \cup V_{i_e+1}$ contains both ends of $e$.
We call each $V_i$ in $(V_1,V_2,...)$ a {\it layer}.
The {\it layered tree-width} of a graph $G$ is the minimum $w$ such that there exist a tree-decomposition of $G$ and a layering of $G$ such that the size of the intersection of any bag and any layer is at most $w$.

We prove the following result, improving a result of Bonamy et al.\ \cite{bbegps} who proved the same for graphs of bounded layered path-width.

\begin{theorem} \label{layered_tw_ad_intro}
For every positive integer $w$, the asymptotic dimension of the class of graphs of layered tree-width at most $w$ is at most 2.
\end{theorem}

Layered tree-width is a common generalization of tree-width and Euler genus of graphs.
A number of classes of graphs with some geometric properties have bounded layered tree-width.
For example, Dujmovi\'{c}, Morin and Wood \cite{dmw} showed that for every nonnegative integer $g$, graphs that can be embedded in a surface of Euler genus at most $g$ have layered tree-width at most $2g+3$.

In fact, classes of graphs of bounded layered tree-width are of interest beyond minor-closed families.
Let $g,k$ be nonnegative integers.
A graph is {\it $(g,k)$-planar} if it can be drawn in a surface of Euler genus at most $g$ such that every edge contains at most $k$ crossings.
So $(g,0)$-planar graphs are exactly the graphs of Euler genus at most $g$.
It is well-known that $(0,1)$-planar graphs (also known as {\it 1-planar graphs} in the literature) can contain an arbitrary graph as a minor.
So the class of $(g,k)$-planar graphs is not a minor-closed family.
On the other hand, Dujmovi\'{c}, Eppstein and Wood \cite{dew} proved that $(g,k)$-planar graphs have layered tree-width at most $(4g+6)(k+1)$.

Hence the following immediate corollary of Theorem \ref{layered_tw_ad_intro} solves a question \cite[Question 4]{bbegps} about $1$-planar graphs and generalizes a result of Fujiwara and Papasoglu \cite{fp} about planar graphs and a result of Bonamy et al. \cite{bbegps} about graphs of bounded Euler genus.

\begin{corollary} \label{gk_planar_ad_intro}
For any nonnegative integers $g$ and $k$, the class of $(g,k)$-planar graphs has asymptotic dimension at most 2.
\end{corollary}

Recall that Corollary \ref{gk_planar_ad_intro} (and hence Theorem \ref{layered_tw_ad_intro}) is optimal since the class of 2-dimensional grids has asymptotic dimension 2.
Other extensively studied graph classes that are known to have bounded layered tree-width include map graphs \cite{dew} and string graphs with bounded maximum degree \cite{djmnw}.
We refer readers to \cite{dew,djmnw} for the discussion for those graphs.

One weakness of layered tree-width is that adding apices can increase layered tree-width a lot.
Note that for any vertex $v$ in a graph $G$ and for any layering of $G$, the neighbors of $v$ must be contained in the union of three consecutive layers.
So if a graph has bounded layered tree-width, then the subgraph induced by the neighbors of any fixed vertex must have bounded tree-width.
However, consider the graphs that can be obtained from 2-dimensional grids by adding a new vertex adjacent to all other vertices.
Since 2-dimensional grids can have arbitrarily large tree-width, such graphs cannot have bounded layered tree-width.

In contrast to the fragility of layered tree-width about adding apices, we show that adding a bounded number of apices does not increase the asymptotic dimension.
Let $\F$ be a class of graphs.
For every nonnegative integer $n$, define $\F^{+n}$ to be the class of graphs such that for every $G \in \F^{+n}$, there exists $Z \subseteq V(G)$ with $\lvert Z \rvert \leq n$ such that $G-Z \in \F$.
We also denote $\F^{+1}$ by $\F^+$. 

\begin{theorem} \label{apex_extension_ad_intro}
Let $\F$ be a class of graphs.
Let $n$ be a nonnegative integer.
Then the asymptotic dimension of $\F^{+n}$ equals the asymptotic dimension of $\F$.
\end{theorem}

This leads to the following strengthening of Theorem \ref{layered_tw_ad_intro}.

\begin{corollary} \label{apex_layered_tw_ad_intro}
Let $k$ be a nonnegative integer.
Let $w$ be a positive integer.
Let $\F$ be a class of graphs such that for every $G \in \F$, there exists $Z \subseteq V(G)$ with $\lvert Z \rvert \leq k$ such that $G-Z$ has layered tree-width at most $w$.
Then the asymptotic dimension of $\F$ is at most 2.
\end{corollary}

The key ingredient of the proof of our aforementioned results in this paper is the following theorem that allows us to show that generating a new class of graphs from hereditary classes by using tree-decompositions of bounded adhesion does not increase the asymptotic dimension.
We say that a class $\F$ of graphs is {\it hereditary} if for every $G \in \F$, every induced subgraph of $G$ belongs to $\F$.
For a graph $G$ and $S \subseteq V(G)$, we denote the subgraph of $G$ induced by $S$ by $G[S]$.

\begin{theorem} \label{tree_extension_intro}
Let $\F$ and $\F'$ be hereditary classes of graphs.
Let $\theta$ be a positive integer.
Let $\F^*$ be a class of graphs such that for every $G \in \F^*$, there exists a tree-decomposition $(T,\X)$ of $G$ of adhesion at most $\theta$, where $\X=(X_t: t \in V(T))$, such that for every $t \in V(T)$, 
	\begin{itemize}
		\item $G[X_t] \in \F$, and 
		\item $\F'$ contains every graph that can be obtained from $G[X_t]$ by adding new vertices and new edges such that for each new vertex $v$, there exists a neighbor $t_v$ of $t$ in $T$ such that the neighbors of $v$ are contained in $X_t \cap X_{t_v}$.
	\end{itemize}
Then $\ad(\F^*) \leq \max\{\ad(\F), \ad(\F'),1\}$. 
\end{theorem}

Theorem \ref{tree_extension_intro} possibly can be further combined with Theorems \ref{minor_ad_intro}, \ref{tw_ad_intro_2}, Corollary \ref{apex_layered_tw_ad_intro}, or other classes with known asymptotic dimension in the literature to generate more complicated classes of graphs without increasing the asymptotic dimension.

The asymptotic dimension is closely related to the notion of weak diameter coloring.
We need some definitions to formally state this relation.

Let $G$ be a graph.
The {\it weak diameter in $G$ of a subset} $S$ of $V(G)$ is the maximum distance in $G$ of two vertices in $S$; the {\it weak diameter in $G$ of a subgraph} $H$ of $G$ is the weak diameter of $V(H)$ in $G$.
Let $k$ be a positive integer.
A {\it $k$-coloring} of $G$ is a function $c:V(G) \rightarrow [k]$.
For a $k$-coloring $c$ of $G$, a {\it $c$-monochromatic component in $G$} is a component of the subgraph of $G$ induced by $\{v \in V(G): c(v)=i\}$ for some $i \in [k]$. 
For an integer $d$, a $k$-coloring $c$ of $G$ is of {\it weak diameter in $G$} at most $d$ if every $c$-monochromatic component in $G$ has weak diameter in $G$ at most $d$.
For every positive integer $\ell$, the {\it $\ell$-th power of $G$}, denoted by $G^\ell$, is the graph obtained from $G$ by adding an edge $xy$ for each pair of distinct vertices of $G$ with distance at most $\ell$.

A simple observation obtained by extending \cite[Observation 1.4]{bbegps} shows a relationship between asymptotic dimension and weak diameter coloring.
Given collections $\U_1,\U_2,...,\U_m$ whose union forms a cover of a metric space $(X,d)$ for defining the asymptotic dimension, we may revise them such that no element of $X$ belongs to members of two distinct $\U_i$'s, so each $\U_i$ gives a color class of an $m$-coloring on $X$, and each member of $\U_i$ is a monochromatic component with bounded weak diameter.

\begin{prop} \label{wd_color_ad_intro}
Let $\F$ be a class of graphs.
Let $m$ be a positive integer.
Then $\ad(\F) \leq m-1$ if and only if there exists a function $f: {\mathbb N} \rightarrow {\mathbb N}$ such that for every $G \in \F$ and $\ell \in {\mathbb N}$, $G^\ell$ is $m$-colorable with weak diameter\footnote{It is worthwhile to mention that for any function $c:V(G) \rightarrow [m]$, $c$ is an $m$-coloring of $G$ and an $m$-coloring of $G^\ell$, but the $c$-monochromatic components in $G$ are different from the $c$-monochromatic components in $G^\ell$.} in $G^\ell$ at most $f(\ell)$.
\end{prop}

An immediate corollary of Proposition \ref{wd_color_ad_intro} is the following.

\begin{corollary} \label{power_ad_intro}
Let $\F$ be a class of graphs.
Let $\ell$ be a positive integer.
Let $\F'$ be the class of graphs such that for every $G \in \F'$, $G=H^\ell$ for some $H \in \F$.
Then $\ad(\F') \leq \ad(\F)$.
\end{corollary}

As pointed out in \cite{bbegps}, weak diameter coloring is also studied under the name ``weak diameter network decomposition'' in distributed computing.
We omit the details about this notion and refer readers to \cite{aglp,bbegps}.

Furthermore, for graphs of bounded maximum degree, weak diameter coloring is equivalent to clustered coloring which is another notion of graph coloring that has attracted wide attention recently (see \cite{w} for a survey).
For positive integers $k$ and $N$, a $k$-coloring of a graph $G$ is said to have {\it clustering $N$} if every $c$-monochromatic component in $G$ contains at most $N$ vertices.
Clearly, a $k$-coloring of clustering $N$ has weak diameter at most $N$.
If the maximum degree of $G$ is at most $\Delta$, and $c$ is a $k$-coloring of $G$ with weak diameter in $G$ at most $w$, then $c$ is a $k$-coloring of clustering $\Delta^{w+1}$.

Hence Theorems \ref{minor_ad_intro} and \ref{tw_ad_intro_2} and Corollary \ref{apex_layered_tw_ad_intro} immediately imply that under the bounded maximum degree condition, $H$-minor free graphs, bounded tree-width graphs, and graphs obtained from bounded layered tree-width graphs by adding a bounded number of apices are 3-colorable, 2-colorable, and 3-colorable with bounded clustering which were proved by the author and Oum \cite{lo}, Alon, Ding, Oporowski and Vertigan \cite{adov}, and the author and Wood \cite{lw}, respectively.
Note that this was observed earlier by Bonamy et al.\ \cite{bbegps}; they pointed out that their \cite{bbegps} results about asymptotic dimension of classes of graphs with bounded maximum degree are enough to deduce the aforementioned results about clustered coloring.

In addition, Theorem \ref{tree_extension_intro} and Corollary \ref{power_ad_intro} give similar results for clustered coloring, as long as the bounded maximum degree condition can be preserved.
For example, as the $\ell$-th power of a graph of maximum degree $\Delta$ has maximum degree at most $\Delta^{\ell+1}$, Theorem \ref{minor_ad_intro} and Corollary \ref{power_ad_intro} imply that for every graph $H$, there exists a constant $k$ such that for every $H$-minor free graph $G$ of maximum degree at most $\Delta$, $G^\ell$ is 3-colorable with clustering $\Delta^{(k+1)(\ell+1)}$. 
Note that the 2-th power of planar graphs can contain any graph as a minor.
So the graph class mentioned above is not minor-closed.
Hence this result is not a special case of \cite{lo}.
After a version of this paper was announced, Wood \cite{w_2} pointed out that the same result can be obtained by combining \cite[Lemma 2]{demww} and \cite[Theorem 17]{dmw_2} but with clustering $O(\ell^{12} \Delta^{5\ell+2} \cdot {\ell+c\Delta \choose c\Delta}^3)$, for some constant $c$ only depending on $H$.

Finally, we remark that our results for asymptotic dimension can be transformed into polynomial time algorithms for finding the covering.
Note that it is equivalent to finding a weak diameter coloring: each $\U_i$ in the covering corresponds to a color class, and each member of $\U_i$ is a monochromatic component.

\begin{theorem} \label{alog_intro}
Let $\F$ be a class of graphs.
Let $w$ and $\ell$ be positive integers.
Let $H$ be a graph.
	\begin{enumerate}
		\item If $\F$ is the class of graphs of tree-width at most $w$, then there exist an integer $N$ and an algorithm such that given $G \in \F$, it finds a 2-coloring of $G^\ell$ with weak diameter at most $N$ in time $O(\lvert V(G) \rvert+\lvert E(G) \rvert)$.
		\item If $\F$ is the class of graphs of layered tree-width at most $w$, then there exist an integer $N$ and an algorithm such that given $G \in \F$ and a pair of layering and tree-decomposition of $G$ witnessing the layered tree-width of $G$, it finds a 3-coloring of $G^\ell$ with weak diameter at most $N$ in time $O(\lvert V(G) \rvert+\lvert E(G) \rvert)$.
		\item If $\F$ is the class of $H$-minor free graphs, then there exist an integer $N$ and an algorithm such that given $G \in \F$, it finds a 3-coloring of $G^\ell$ with weak diameter at most $N$ in time $O(\lvert V(G) \rvert^2)$.
	\end{enumerate}
\end{theorem}

\section{Proof sketch and organization of the paper}

In this section we sketch our proofs of the results mentioned in Section \ref{sec:intro} and explain the organization of this paper.
As discussed in the previous section, it suffices to prove Theorems \ref{minor_ad_intro}, \ref{tw_ad_intro_2}, \ref{layered_tw_ad_intro}, \ref{apex_extension_ad_intro}, \ref{tree_extension_intro}, \ref{alog_intro} and Proposition \ref{wd_color_ad_intro}.

Proposition \ref{wd_color_ad_intro} easily follows from the definition of the asymptotic dimension, but it is useful as it allows us to use the terminology for coloring.
From now on, it suffices to prove that given a class $\F$ and a positive integer $\ell$, there exists $N$ such that for every $G \in \F$, $G^\ell$ is $m$-colorable with weak diameter in $G^\ell$ at most $N$, where $m \in \{2,3\}$ is the number of colors required for the theorems.
We include a detailed proof of Proposition \ref{wd_color_ad_intro} in Section \ref{sec:ad_wd_color}.

\subsection{Proof of the existence of coloring}
In this subsection we sketch proofs of Theorems \ref{minor_ad_intro}, \ref{tw_ad_intro_2}, \ref{layered_tw_ad_intro}, \ref{apex_extension_ad_intro} and \ref{tree_extension_intro}.
Theorem \ref{alog_intro} can be proved by transforming the arguments in this subsection into algorithms, and we will discuss it in Section \ref{subsec:alog}.

\subsubsection{Machinery for extending precoloring}

In order to prove that $G^\ell$ is $m$-colorable with bounded weak diameter in this paper, we usually prove a stronger statement which roughly says that we can always extend a precoloring on a ``reasonable'' set to a coloring with bounded weak diameter of the entire graph.
Here ``reasonable'' sets are the sets that are not far from a set of bounded size.

We need some definitions to make it formal.
Let $r$ and $k$ be nonnegative integers.
Let $G$ be a graph.
For $S \subseteq V(G)$, we define $N_G^{\leq r}[S] = \{v \in V(G):$ there exists a path in $G$ from $v$ to $S$ with length\footnote{The {\it length} of a path is the number its edges.} at most $r\}$.
We say that $Z \subseteq V(G)$ is {\it $(k,r)$-centered} if there exists $S \subseteq V(G)$ with $\lvert S \rvert \leq k$ such that $Z \subseteq N_G^{\leq r}[S]$.

The first technical tool in this paper is Lemma \ref{deleting_centered_set} stating that if $Z$ is a $(k,r)$-centered set in $G$, then combining any $m$-coloring $c_Z$ on $Z$ and any $m$-coloring $c$ of $(G-Z)^\ell$ with bounded weak diameter in $(G-Z)^\ell$ gives an $m$-coloring of $G^\ell$ with bounded weak diameter in $G^\ell$.
The proof of this lemma follows from induction on $k$.
Let $S$ be a set with $\lvert S \rvert \leq k$ such that $Z \subseteq N_G^{\leq r}[S]$.
When there exist two vertices in $S$ not very far from each other in $G^\ell$, we can remove some vertex from $S$ to show that $Z$ is a $(k-1,r')$-centered set for some $r'>r$, so we are done by induction.
So we may assume that the vertices in $S$ are pairwise far apart in $G^\ell$.
Hence every $c$-monochromatic component $M$ of $(G-Z)^\ell$ can be close in $G$ to at most one vertex in $S$, as $c$ has small weak diameter.
If $M$ is not close to any vertex in $S$, then $M$ is far from $Z$, so $M$ remains a monochromatic component in $c_Z \cup c$.
If $M$ is close to one vertex $s$ in $S$, then it can be enlarged into a $(c_Z \cup c)$-monochromatic component $M'$, but $M'$ is still not far from $s$, so the weak diameter remains small.

So we complete the sketch of the proof of Lemma \ref{deleting_centered_set}.
Section \ref{sec:center} is dedicated to proving Lemma \ref{deleting_centered_set} and its simple applications.

\subsubsection{Using Theorem \ref{tree_extension_intro}}

The other technical tool in this paper is Theorem \ref{tree_extension_intro}.
Before sketching the proof of Theorem \ref{tree_extension_intro}, we explain how to derive all remaining results by combining Lemma \ref{deleting_centered_set} and Theorem \ref{tree_extension_intro}.
These applications of Lemma \ref{deleting_centered_set} and Theorem \ref{tree_extension_intro} are included in Section \ref{sec:app_ad}.

	\begin{itemize}
		\item Theorem \ref{apex_extension_ad_intro} immediately follows from Lemma \ref{deleting_centered_set}, since every graph in $\F^{+n}$ can be made a graph in $\F$ by deleting a set of size at most $n$ which is an $(n,0)$-centered set.
		\item Theorem \ref{apex_extension_ad_intro} immediately implies that classes of graphs with a vertex-cover of bounded size have asymptotic dimension 0, since this class is $\F^{+n}$ for some $n$, where $\F$ is the set of edgeless graphs.
		\item Then the bounded tree-width result (Theorem \ref{tw_ad_intro_2}) immediately follows from Theorem \ref{tree_extension_intro}.
			If we set the class $\F^*$ in Theorem \ref{tree_extension_intro} to be the class of graphs of bounded tree-width, then the corresponding classes $\F$ and $\F'$ mentioned in Theorem \ref{tree_extension_intro} are the class of graphs with a bounded number of vertices and the class of graphs with a vertex-cover of bounded size, respectively.
		\item Using \cite[Theorem 5.2]{bbegps} (Theorem \ref{layer_trick} in this paper), the bounded layered tree-width result (Theorem \ref{layered_tw_ad_intro}) immediately follows from Theorem \ref{tw_ad_intro_2}.
		\item Then Corollary \ref{apex_layered_tw_ad_intro} immediately follows from Theorems \ref{layered_tw_ad_intro} and \ref{apex_extension_ad_intro}.
		\item Structure theorems for excluding minors developed in \cite{rs_XVI} and \cite{dmw} state that every $H$-minor free graph has a tree-decomposition of bounded adhesion such that every torso can be made a graph of bounded layered tree-width by deleting a bounded number of vertices.
			If we set the class $\F^*$ in Theorem \ref{tree_extension_intro} to be the class of $H$-minor free graphs, then the corresponding $\F$ in Theorem \ref{tree_extension_intro} is a class mentioned in Corollary \ref{apex_layered_tw_ad_intro}.
			It can be shown that the corresponding $\F'$ in Theorem \ref{tree_extension_intro} is also a class mentioned in Corollary \ref{apex_layered_tw_ad_intro}. (See the proof of Theorem \ref{minor_ad} for details.)
			So Theorem \ref{minor_ad_intro} follows from Theorem \ref{tree_extension_intro} and Corollary \ref{apex_layered_tw_ad_intro}.
	\end{itemize}

\subsubsection{Proof sketch of Theorem \ref{tree_extension_intro}}

Therefore, it suffices to prove Theorem \ref{tree_extension_intro}, and it is the most technical part of the paper.
Section \ref{sec:tree} is devoted to this task. 

Let $G \in \F^*$ and $(T,\X)$ be a tree-decomposition of $G$ as stated in Theorem \ref{tree_extension_intro}.
As $(T,\X)$ has bounded adhesion, we can treat $T$ as a rooted tree, and the bag of the root, denoted by $X_{t^*}$, has size at most the adhesion by creating a redundant bag if necessary.
We shall prove a stronger statement: for every $Z \subseteq N_G^{\leq 3\ell}[X_{t^*}]$, every precoloring $c$ on $Z$ with at most $m$ colors extends to an $m$-coloring of $G^\ell$ with bounded weak diameter in $G^\ell$, by induction on the adhesion of $(T,\X)$, and subject to this, induction on $\lvert V(G) \rvert$. (See Lemma \ref{tree_extension} for a precise statement.)

By first extending $c$ to $N_G^{\leq 3\ell}[X_{t^*}]$, we may assume $Z=N_G^{\leq 3\ell}[X_{t^*}]$.
Hence the subgraph $T_0$ of $T$ induced by the nodes whose bags intersect $Z$ is a subtree of $T$ containing $t^*$.
Let $G_0$ be the subgraph of $G$ induced by the bags of the nodes in $T_0$.
Let $U_E$ be the set of edges of $T$ with exactly one end in $T_0$.
For every $e \in U_E$, let $G_e$ be the subgraph of $G$ induced by the bags of the nodes in the component of $T-e$ disjoint from $T_0$.
Then for each $e \in U_E$, we can partition $V(G_0) \cap V(G_e)$ into sets such that any two vertices in $V(G_0) \cap V(G_e)$ are not far from each other in $G_e$ if and only if they are contained in the same part in this partition, denoted by $\P_e$.
Note that it can be done as $\lvert V(G_0) \cap V(G_e) \rvert$ is bounded by the adhesion of $(T,\X)$.

Note that $Z$ is an $(\lvert X_{t^*} \rvert,3\ell)$-centered set, and $G_0-Z$ has a tree-decomposition of smaller adhesion by the definition of $G_0$.
So the precoloring $c$ on $Z$ restricted to $Z-Z=\emptyset$ can be extended to $(G_0-Z)^\ell$ by induction.
Hence the precoloring $c$ on $Z$ can be extended to $G_0^\ell$ by Lemma \ref{deleting_centered_set}.
However, it is troublesome to further extend the coloring to $G^\ell$, as no information about $G_e$ for $e \in U_E$ can be seen from $G_0$ and edges of $G^\ell$ with ends in $V(G_0)$ cannot be completely told from $G_0^\ell$.
To overcome this difficulty, we add ``gadgets'' to $G_0$ to obtain a graph $G_0'$ such that extending $c$ from $Z$ to $(G_0')^\ell$ gives sufficient information about how to further extend it to $G^\ell$.
The gadgets we added to form $G_0'$ are a vertex $v_Y$ for each $e \in U_E$ and each part $Y \in \P_e$, and the edges between $v_Y$ and $Y$.

However, $G_0'-Z$ possibly does not have a tree-decomposition of smaller adhesion, so the induction hypothesis cannot be applied to $G_0'-Z$.
Instead, we setup a more technical induction hypothesis to overcome this difficulty.
This is the motivation of $(\eta,\theta,\F,\F')$-constructions mentioned in Section \ref{sec:tree}.
So we can extend $c$ to an $m$-coloring of $(G_0')^\ell$ by using this technical setting.

Note that no vertex in $Z$ is in $\bigcup_{e \in U_E}G_e$.
Then for each $e \in U_E$, we color vertices in $G_e-V(G_0)$ that have distance in $G_e$ at most $\ell$ from $V(G_e) \cap V(G_0)$ according to the colors on $v_Y$ for $Y \in \P_e$.
Call the set of vertices colored in this step $Z_1$.
Then we color every uncolored vertex in $G_e$ that have distance in $G_e$ at most $\ell$ from $Z_1$ color 1.
Call the set of vertices colored in this step $Z_2$.
Then we color every uncolored vertex in $G_e$ that have distance in $G_e$ at most $\ell$ from $Z_2$ color 2.
It makes sure that no matter how we further color other vertices, every monochromatic component intersecting $G_0$ must be contained in $V(G_0) \cup Z_1 \cup Z_2$; and the information of the coloring on $v_Y$ helps us show that such monochromatic components have small weak diameter.

At this point, for every $e \in U_E$, the vertices colored in $G_e$ are contained in $N_G^{\leq 3\ell}[V(G_0) \cap V(G_e)]$.
Hence for each $e \in U_E$, we can extend this precoloring to an $m$-coloring of $G_e^\ell$ with bounded weak diameter in $G_e^\ell$ by induction, since $G_e$ has fewer vertices than $G$.
Every monochromatic component not intersecting $G_0$ must be contained in $G_e$ for some $e \in U_E$ and hence has small weak diameter.
This completes the sketch of the proof of Theorem \ref{tree_extension_intro} (or Lemma \ref{tree_extension}).

\subsection{Algorithmic aspects} \label{subsec:alog}

In this subsection we prove Theorem \ref{alog_intro} by transforming the arguments in the previous subsection into algorithms.

\subsubsection{Algorithm for Lemma \ref{deleting_centered_set}}

Lemma \ref{deleting_centered_set} shows that simply combining an $m$-coloring of $Z$ and an $m$-coloring of $(G-Z)^\ell$ gives an $m$-coloring of $G^\ell$ with bounded weak diameter, so a desired $m$-coloring of $G^\ell$ can be obtained in linear time as long as the colorings of $Z$ and $(G-Z)^\ell$ are given.

\subsubsection{Time complexity for Theorem \ref{tree_extension_intro}}

Now we consider the time complexity for applying Theorem \ref{tree_extension_intro}.
Assume that there exists a function $h$ such that for every graph $G \in \F \cup \F'$, an $m$-coloring of $G^\ell$ with bounded weak diameter in $G^\ell$ can be found in time $h(\lvert V(G) \rvert + \lvert E(G) \rvert)$.
Our proof of Theorem \ref{tree_extension_intro} (or Lemma \ref{tree_extension}) implies that if a tree-decomposition $(T,\X)$ of a graph $G \in \F^*$ witnessing the membership of $G$ in $\F^*$ (or more precisely, an $(\eta,\theta,\F,\F')$-construction of $G$ mentioned in Section \ref{sec:tree}) is given, then an $m$-coloring of $G^\ell$ with bounded weak diameter in $G^\ell$ can be found in time $O(h(2^{\theta^2} \cdot (\lvert V(G) \rvert + \lvert E(G) \rvert)))$.

\subsubsection{Proof of Theorem \ref{alog_intro}}

First, we show the time complexity for our algorithm for the bounded tree-width classes (Statement 1 of Theorem \ref{alog_intro}).
Let $\F^*$ be the class of graphs of tree-width at most $w$.
The corresponding classes $\F$ and $\F'$ in Theorem \ref{tree_extension_intro} are the class of graphs on at most $w+1$ vertices and the class of graphs that have a vertex-cover of size at most $w+1$.
Hence $\F \subseteq \F'$.
Given a graph $L \in \F'$, we can obtain a vertex-cover $S$ of $L$ of size at most $2(w+1)$ in linear time by collecting the ends of the edges in a maximal matching.
Then by coloring every vertex in $S$ color 1 and coloring every vertex in $L-S$ color 2, we obtain a 2-coloring of $L$ with bounded weak diameter by Lemma \ref{deleting_centered_set}, since $S$ is a $(2w+2,0)$-centered set.
Hence the function $h$ is a linear function.
By a result of Bodlaender \cite{b}, given a graph $G \in \F^*$, one can find a tree-decomposition of $G$ of width at most $w$ in linear time.
Since $h$ is a linear function, the algorithm given by the proof of Theorem \ref{tree_extension_intro} shows that one can find a 2-coloring of $G^\ell$ with bounded weak diameter in linear time.
This proves Statement 1 of Theorem \ref{alog_intro}.

Next, we show the time complexity for the bounded layered tree-width case (Statement 2 of Theorem \ref{alog_intro}).
Assume that $G$ is a graph of layered tree-width at most $w$, and a pair of witnessing layering and tree-decomposition is given.
Since a 2-coloring with bounded weak diameter of the $\ell'$-th power of a graph of tree-width at most $w$ can be found in linear time for any fixed $\ell'$, by the machinery developed in \cite{bbegps} (or \cite{bdlm}), a 3-coloring with bounded weak diameter of $G^\ell$ can be found in linear time.
This proves Statement 2 of Theorem \ref{alog_intro}.

Finally, we show the time complexity for the class of $H$-minor free graphs (Statement 3 of Theorem \ref{alog_intro}).
Let $\F^*$ be the class of graphs of $H$-minor free graphs.
By a result of Grohe, Kawarabayashi and Reed \cite{gkr}, there exist a quadratic time algorithm and an integer $p$ (only depending on $H$) such that given a graph $G \in \F^*$, one can find a tree-decomposition $(T,\X)$ of $G$ of adhesion at most $p$ such that for every $t \in V(T)$, one can find $Z_t \subseteq V(G_t)$ with $\lvert Z_t \rvert \leq p$ and a ``$p$-nearly embedding'' of $G_t-Z_t$, where $G_t$ is the graph obtained from $G[X_t]$ by adding edges such that $X_t \cap X_{t'}$ is a clique for each neighbor $t'$ of $t$ in $T$. 
Since a $p$-nearly embedding of $G_t-Z_t$ is given, the proof of a result of Dujmovi\'{c}, Morin and Wood \cite{dmw} shows that there exists an integer $q$ (only depending on $p$) such that for every $t \in V(T)$, $G_t-Z_t$ has layered tree-width at most $q$, and a pair of witnessing layering and tree-decomposition of $G_t-Z_t$ can be found in quadratic time.
So the corresponding class $\F$ in Theorem \ref{tree_extension_intro} is $\W_q^{+p}$, where $\W_q$ is the class of layered tree-width at most $q$.
Lemma \ref{small_exten_minor} shows that the corresponding $\F'$ in Theorem \ref{tree_extension_intro} is $\W_{q+1}^{+p}$, and whenever a pair of witnessing layering and tree-decomposition of $G_t-Z_t$ is given, one can find in linear time a set $Z_t' \subseteq V(G_t')$ with $\lvert Z_t' \rvert \leq p$ and a pair of layering and tree-decomposition of $G_t'-Z_t'$ witnessing the membership for $G_t'$ in $\W_{q+1}^{+p}$, where $G_t'$ is the corresponding member of $\F'$ generated by $G[X_t]$ stated in the statement of Theorem \ref{tree_extension_intro}.
By Lemma \ref{deleting_centered_set} and Statement 2 of Theorem \ref{alog_intro}, the function $h$ for the time complexity is quadratic in the number of edges.
Note that every graph in $\W_{q+1}^{+p}$ is $(3q+3+p)$-degenerate, so $h$ is quadratic in the number of vertices.
Therefore, there exists a quadratic time algorithm to find a 3-coloring with bounded weak diameter of the $\ell$-th power of an $H$-minor free graph.
This proves Statement 3 of Theorem \ref{alog_intro}.

\section{Asymptotic dimension and weak diameter coloring} \label{sec:ad_wd_color}

\begin{prop} \label{wd_color_ad}
Let $\F$ be a class of graphs.
Let $m$ be a positive integer.
Then the following are equivalent.
	\begin{enumerate}
		\item $\ad(\F) \leq m-1$.
		\item There exists a function $f: {\mathbb R}^+ \rightarrow {\mathbb R}^+$ such that for every $G \in \F$ and $r \in {\mathbb R}^+$, there exist $m$ collections $\X_1,\X_2,...,\X_m$ such that 
			\begin{itemize}
				\item $\bigcup_{i=1}^m\bigcup_{X \in \X_i}X \supseteq V(G)$,
				\item for any $i \in [m]$ and $X \in \X_i$, the weak diameter of $X$ in $G$ is at most $f(r)$, and
				\item for any $i \in [m]$, distinct $X,X' \in \X_i$ and elements $x \in X$ and $x' \in X'$, the distance between $x$ and $x'$ in $G$ is greater than $r$.
			\end{itemize}
		\item There exists a function $f: {\mathbb N} \rightarrow {\mathbb N}$ such that for every $G \in \F$ and $\ell \in {\mathbb N}$, there exist $m$ collections $\X_1,\X_2,...,\X_m$ such that 
			\begin{itemize}
				\item $\bigcup_{i=1}^m\bigcup_{X \in \X_i}X \supseteq V(G)$,
				\item for any $i \in [m]$ and $X \in \X_i$, the weak diameter of $X$ in $G$ is at most $f(\ell)$, and
				\item for any $i \in [m]$, distinct $X,X' \in \X_i$ and elements $x \in X$ and $x' \in X'$, the distance between $x$ and $x'$ in $G$ is greater than $\ell$.
			\end{itemize}
		\item There exists a function $f: {\mathbb N} \rightarrow {\mathbb N}$ such that for every $G \in \F$ and $\ell \in {\mathbb N}$, $G^\ell$ is $m$-colorable with weak diameter in $G$ at most $f(\ell)$.
		\item There exists a function $f: {\mathbb N} \rightarrow {\mathbb N}$ such that for every $G \in \F$ and $\ell \in {\mathbb N}$, $G^\ell$ is $m$-colorable with weak diameter in $G^\ell$ at most $f(\ell)$.
	\end{enumerate}
\end{prop}

\begin{pf}
Statements 1 and 2 are equivalent by the definition of the asymptotic dimension.
Statements 2 and 3 are equivalent since the distance between any two vertices in a graph is an integer.

Now we show that Statement 3 implies Statement 4.
Fix $G$ and $\ell$.
For every $v \in V(G)$, define $c(v)$ to be the minimum $i \in [m]$ such that $v \in \bigcup_{X \in \X_i}X$.
Then $c$ is an $m$-coloring of $G$ and hence an $m$-coloring of $G^\ell$.
Let $M$ be a $c$-monochromatic component in $G^\ell$.
For every edge $uv \in E(M)$, since $M \subseteq G^\ell$, the distance between $u,v$ in $G$ is at most $\ell$, so there exists $X \in \X_{c(u)}$ such that $X$ contains $u$ and $v$.
Hence there exists $X \in \X_{c(M)}$ such that $V(M) \subseteq X$.
Since the weak diameter of $X$ in $G$ is at most $f(\ell)$, so is $V(M)$.
Hence $G^\ell$ is $m$-colorable with weak diameter in $G$ at most $f(\ell)$.
So Statement 3 implies Statement 4.

Now we show that Statement 4 implies Statement 3.
Fix $G$ and $\ell$.
Let $c$ be an $m$-coloring of $G^\ell$ with weak diameter in $G$ at most $f(\ell)$.
For each $i \in [m]$, define $\X_i=\{V(M): M$ is a $c$-monochromatic component in $G^\ell$ with $c(M)=i\}$.
So $\bigcup_{i=1}^m\bigcup_{X \in \X_i}X \supseteq V(G)$.
For each $i \in [m]$ and $X \in \X_i$, the weak diameter of $X$ in $G$ is the weak diameter in $G$ of some $c$-monochromatic component in $G^\ell$, so it is at most $f(\ell)$.
For each $i \in [m]$, distinct $X,X' \in \X_i$ and elements $x \in X$ and $x' \in X'$, $x$ and $x'$ belong to different $c$-monochromatic components in $G^\ell$ with the same color, so they are not adjacent in $G^\ell$, and hence the distance between $x$ and $x'$ in $G$ is greater than $\ell$.
So Statement 4 implies Statement 3.

Now we show that Statements 4 and 5 are equivalent.
Note that a set $S$ in a graph $G$ has weak diameter $k$ in $G$ implies $S$ has weak diameter at most $k$ in $G^\ell$ since $G \subseteq G^\ell$ for any positive integer $\ell$.
A set $S$ in a graph $G$ has weak diameter $k$ in $G^\ell$ implies that for any two vertices in $S$, there exists a path in $G^\ell$ of length at most $k$ between them, so there exists a walk in $G$ of length at most $\ell k$ between them; hence the weak diameter of $S$ in $G$ is at most $\ell k$.
Therefore Statements 4 and 5 are equivalent.
\end{pf}

\section{Centered sets} \label{sec:center}

Let $r$ and $k$ be nonnegative integers.
Let $G$ be a graph.
Let $S \subseteq V(G)$.
Recall that we define $N_G^{\leq r}[S] = \{v \in V(G):$ there exists a path in $G$ from $v$ to $S$ with length at most $r\}$.

For $i \in [2]$, let $f_i$ be a function with domain $S_i$.
If $f_1(x)=f_2(x)$ for every $x \in S_1 \cap S_2$, then we define $f_1 \cup f_2$ to be the function with domain $S_1 \cup S_2$ such that for every $x \in S_1 \cup S_2$, $(f_1 \cup f_2)(x) = f_{i_x}(x)$, where $i_x$ is an element in $[2]$ such that $x \in S_{i_x}$.

\begin{lemma} \label{deleting_centered_set}
For any nonnegative integers $k,r$ and positive integers $\ell,N$, there exists an integer $N^* \geq N$ such that the following holds.
Let $G$ be a graph.
Let $S \subseteq V(G)$ with $\lvert S \rvert \leq k$.
Let $Z \subseteq N_G^{\leq r}[S]$.
Let $m$ be a positive integer.
Let $c_Z: Z \rightarrow [m]$.
Let $c$ be an $m$-coloring of $(G-Z)^\ell$ with weak diameter in $(G-Z)^\ell$ at most $N$.
Then the $m$-coloring $c \cup c_Z$ of $G^\ell$ has weak diameter in $G^\ell$ at most $N^*$.
\end{lemma}

\begin{pf}
 We define $f: ({\mathbb N} \cup \{0\}) \times ({\mathbb N} \cup \{0\}) \times {\mathbb N} \times {\mathbb N} \rightarrow {\mathbb N}$ to be the function such that 
 	\begin{itemize}
		\item $f(0,x,y,z)= z$, and 				
                \item for every $\alpha \in {\mathbb N}$, $f(\alpha,x,y,z) = 2x+2y+2f(\alpha-1,x,y,z)$. 
 	\end{itemize}

It suffices to show that the $m$-coloring $c \cup c_Z$ of $G^\ell$ has weak diameter in $G^\ell$ at most $f(k,r,\ell,N)$. 
We shall prove it by induction on $k$.
When $k=0$, $S=Z=\emptyset$, and hence we are done since $f(0,r,\ell,N) \geq N$.
So we may assume that $k \geq 1$ and this lemma holds when $k$ is smaller.

Let $s_1 \in S$.
Let $S'=S-\{s_1\}$.
Let $Z_1 = Z \cap N_G^{\leq r}[s_1] - N_G^{\leq r}[S']$.
Let $Z' = Z-Z_1$.
So $Z' \subseteq N_G^{\leq r}[S']$.
If $Z' \not \subseteq N_{G-Z_1}^{\leq r}[S']$, then there exists a path $Q$ in $G$ from an element $s' \in S'$ to a vertex $Z'$ with length in $G$ at most $r$ such that $Z_1 \cap V(Q) \neq \emptyset$, so $N_G^{\leq r}[S'] \cap Z_1 \supseteq V(Q) \cap Z_1 \neq \emptyset$, a contradiction. 
So $Z' \subseteq N_{G-Z_1}^{\leq r}[S']$.
Since $(G-Z)^\ell = ((G-Z_1)-Z')^\ell$, $c$ is an $m$-coloring of $((G-Z_1)-Z')^\ell$ with weak diameter in $((G-Z_1)-Z')^\ell$ at most $N$.
Since $\lvert S' \rvert < \lvert S \rvert$, by the induction hypothesis, the $m$-coloring $c \cup c_Z|_{Z'}$ of $(G-Z_1)^\ell$ has weak diameter in $(G-Z_1)^\ell$ at most $f(k-1,r,\ell,N)$.

Let $c'=c \cup c_Z|_{Z'}$.
Note that $c \cup c_Z = c' \cup c_Z|_{Z_1}$.

Suppose to the contrary that the weak diameter in $G^\ell$ of $c \cup c_Z$ is greater than $f(k,r,\ell,N)$.
Then there exists a $(c \cup c_Z)$-monochromatic component $M$ in $G^\ell$ with weak diameter in $G^\ell$ greater than $f(k,r,\ell,N)$.
Hence there exist vertices $u$ and $v$ of $M$ such that the distance in $G^\ell$ between $u$ and $v$ is greater than $f(k,r,\ell,N)$.

Since $\ell \geq 1$, $G \subseteq G^\ell$.
If both $u$ and $v$ are contained in $Z_1$, then the distance in $G$ between $u$ and $v$ is at most $2r$, so the distance in $G^\ell$ between $u$ and $v$ is at most $2r \leq f(k,r,\ell,N)$, a contradiction.
So at least one of $u$ and $v$, say $u$, is not in $Z_1$.
Hence there exists a $c'$-monochromatic component $M_u$ in $(G-Z_1)^\ell$ containing $u$.
Note that $M_u \subseteq M$.

\noindent{\bf Claim 1:} The distance in $G^\ell$ between $s_1$ and $u$ is at most $r+\ell+f(k-1,r,\ell,N)$.

\noindent{\bf Proof of Claim 1:}
Suppose to the contrary that the distance in $G^\ell$ between $s_1$ and $u$ is greater than $r+\ell+f(k-1,r,\ell,N)$.

Suppose to the contrary that there exists a vertex $z \in Z_1$ adjacent in $G^\ell$ to a vertex $u'$ in $V(M_u)$.
So the distance in $G^\ell$ between $z$ and $u'$ is at most $\ell$.
Since $c'$ has weak diameter in $(G-Z_1)^\ell$ at most $f(k-1,r,\ell,N)$ and $M_u$ contains $u$ and $u'$, the distance in $(G-Z_1)^\ell$ between $u$ and $u'$ is at most $f(k-1,r,\ell,N)$.
Since $(G-Z_1)^\ell \subseteq G^\ell$, the distance in $G^\ell$ between $u$ and $u'$ is at most $f(k-1,r,\ell,N)$.
Since the distance in $G^\ell$ between $z$ and $u'$ is at most $\ell$, the distance in $G^\ell$ between $s_1$ and $u$ is at most $r+\ell+f(k-1,r,\ell,N)$, a contradiction.

Hence there exists no vertex in $Z_1$ adjacent in $G^\ell$ to  $V(M_u)$.
In particular, there exists no edge of $M-E(M_u)$ incident with $V(M_u)$.
Hence, $M=M_u$.
So $u$ and $v$ are contained in the same $c'$-monochromatic component.
Since $(G-Z_1)^\ell \subseteq G^\ell$, the distance in $G^\ell$ between $u$ and $v$ is at most $f(k-1,r,\ell,N) \leq f(k,r,\ell,N)$, a contradiction.
$\Box$

Therefore, the distance in $G^\ell$ between $s_1$ and $u$ is at most $r+\ell+f(k-1,r,\ell,N)$ by Claim 1.
Similarly, either $v \in Z_1$, or the distance in $G^\ell$ between $s_1$ and $v$ is at most $r+\ell+f(k-1,r,\ell,N)$.
If $v \in Z_1$, then since $G \subseteq G^\ell$, $v \in N_G^{\leq r}[\{s_1\}] \subseteq N_{G^\ell}^{\leq r}[\{s_1\}]$.
So the distance in $G^\ell$ between $s_1$ and $v$ is at most $r+\ell+f(k-1,r,\ell,N)$ in either case.
Therefore, the distance in $G^\ell$ between $u$ and $v$ is at most $2r+2\ell+2f(k-1,r,\ell,N) \leq f(k,r,\ell,N)$, a contradiction.
This proves the lemma.
\end{pf}

\bigskip

Let $r$ and $k$ be nonnegative integers.
Let $G$ be a graph.
Let $Z \subseteq V(G)$.
Recall that we say that $Z$ is {\it $(k,r)$-centered} if there exists $S \subseteq V(G)$ with $\lvert S \rvert \leq k$ such that $Z \subseteq N_G^{\leq r}[S]$.

Let $\ell,N,m$ be positive integers.
We say a class $\F$ of graphs is {\it $(m,\ell,N)$-nice} if for every $G \in \F$, $G^\ell$ is $m$-colorable with weak diameter in $G^\ell$ at most $N$.

\begin{lemma} \label{apex_extension}
For any positive integers $\ell,N$ and nonnegative integer $n$, there exists $N^*$ such that the following holds.
Let $m$ be a positive integer.
If $\F$ is an $(m,\ell,N)$-nice class of graphs, then $\F^{+n}$ is an $(m,\ell,N^*)$-nice class. 
\end{lemma}

\begin{pf}
Let $\ell,N$ be positive integers and let $n$ be an nonnegative integer.
Define $N^*$ to be the integer $N^*$ mentioned in Lemma \ref{deleting_centered_set} by taking $(k,r,\ell,N)=(n,0,\ell,N)$.

Let $G \in \F^{+n}$.
So there exists $Z \subseteq V(G)$ with $\lvert Z \rvert \leq n$ such that $G-Z \in \F$.
Since $\F$ is an $(m,\ell,N)$-nice class, there exists an $m$-coloring $c$ of $(G-Z)^\ell$ with weak diameter in $(G-Z)^\ell$ at most $N$. 
Let $c_Z$ be an $m$-coloring of $Z$.
Since $Z$ is $(n,0)$-centered, by Lemma \ref{deleting_centered_set}, $(c \cup c_Z)$ is an $m$-coloring of $G^\ell$ with weak diameter in $G^\ell$ at most $N^*$.

Therefore, $\F^{+n}$ is $(m,\ell,N^*)$-nice.
\end{pf}

\bigskip

Now we are ready to prove Theorem \ref{apex_extension_ad_intro}.
The following is a restatement.

\begin{theorem} \label{apex_extension_ad}
Let $\F$ be a class of graphs. 
Let $n$ be a nonnegative integer.
Then $\ad(\F^{+n})=\ad(\F)$.
\end{theorem}

\begin{pf}
By Proposition \ref{wd_color_ad}, there exists a function $f: {\mathbb N} \rightarrow {\mathbb N}$ such that for every $\ell \in {\mathbb N}$, $\F$ is $(\ad(\F)+1,\ell,f(\ell))$-nice.
By Lemma \ref{apex_extension}, there exists a function $g: {\mathbb N} \rightarrow {\mathbb N}$ such that for every $\ell \in {\mathbb N}$, $\F^{+n}$ is $(\ad(\F)+1,\ell,g(f(\ell)))$-nice.
By Proposition \ref{wd_color_ad}, $\ad(\F^{+n}) \leq \ad(\F)$.
Since $\F^{+n} \supseteq \F$, $\ad(\F^{+n})=\ad(\F)$.
\end{pf}

\begin{lemma} \label{all_centered}
For any nonnegative integers $k$ and $r$ and positive integer $\ell$, there exists a positive integer $N$ such that the following holds.
Let $G$ be a graph.
Let $m$ be a positive integer.
If $V(G)$ is $(k,r)$-centered, then any $m$-coloring of $G^\ell$ has weak diameter in $G^\ell$ at most $N$. 
\end{lemma}

\begin{pf}
Let $k,r$ be nonnegative integers, and let $\ell$ be a positive integer.
Define $N$ to be the integer $N^*$ mentioned in Lemma \ref{deleting_centered_set} by taking $(k,r,\ell,N)=(k,r,\ell,1)$ in Lemma \ref{deleting_centered_set}.
Since $G-V(G)$ has no vertex, every $m$-coloring of $G-V(G)$ has weak diameter at most 0.
Since $V(G)$ is $(k,r)$-centered, for every function $c: V(G) \rightarrow [m]$, $c$ is an $m$-coloring of $G^\ell$ with weak diameter in $G^\ell$ at most $N$ by Lemma \ref{deleting_centered_set}.
\end{pf}

\bigskip

A {\it vertex-cover} of a graph $G$ is a subset $S$ of $V(G)$ such that $V(G)-S$ has no edge.

\begin{lemma} \label{vc_color}
For any nonnegative integer $k$ and positive integer $\ell$, there exists a positive integer $N$ such that the following holds.
Let $m$ be a positive integer.
If $G$ is a graph that has a vertex-cover of size at most $k$, then any $m$-coloring of $G^\ell$ has weak diameter in $G^\ell$ at most $N$.
\end{lemma}

\begin{pf}
Let $k$ be a nonnegative integer, and let $\ell$ be a positive integer.
Define $N$ to be the positive integer $N$ mentioned in Lemma \ref{all_centered} by taking $(k,r,\ell)=(k,1,\ell)$.

Let $m$ be a positive integer.
Let $G$ be a graph that has a vertex-cover $S$ of size at most $k$.
Let $G_1=G[N_G^{\leq 1}[S]]$.
Since $V(G_1)$ is $(k,1)$-centered, any $m$-coloring of $G_1^\ell$ has weak diameter in $G_1^\ell$ at most $N$ by Lemma \ref{all_centered}.
Since $S$ is a vertex-cover, every component of $G^\ell$ either is contained in $G_1^\ell$ or consists of one vertex.
So any $m$-coloring of $G^\ell$ has weak diameter in $G^\ell$ at most $N$.
\end{pf}

\begin{lemma} \label{vc_ad}
Let $k$ be a nonnegative integer.
Let $\F$ be a class of graphs such that every graph in $\F$ has a vertex-cover of size at most $k$.
Then $\ad(\F) \leq 0$.
\end{lemma}

\begin{pf}
Let $k$ be a nonnegative integer.
Define $f: {\mathbb N} \rightarrow {\mathbb N}$ to be the function such that for every $x \in {\mathbb N}$, $f(x)$ equals the integer $N$ mentioned in Lemma \ref{vc_color} by taking $(k,\ell)=(k,x)$.
Then for every graph $G \in \F$, since $G$ has a vertex-cover of size at most $k$, by Lemma \ref{vc_color}, for every $\ell \in {\mathbb N}$, $G^\ell$ is $1$-colorable with diameter in $G^\ell$ at most $f(\ell)$.
By Proposition \ref{wd_color_ad}, $\ad(\F) \leq 0$.
\end{pf}

\section{Gluing along a tree} \label{sec:tree}

Let $G$ be a graph, and let $(T,\X)$ be a tree-decomposition of $G$, where $\X=(X_t: t \in V(T))$.
For every $S \subseteq V(T)$, we define $X_S = \bigcup_{t \in S}X_t$; for every subgraph $H$ of $T$, we define $X_H = X_{V(H)}$.

A {\it rooted tree} is a directed graph whose underlying graph is a tree such that there exists a unique vertex $t$ with in-degree 0.
We call $t$ the {\it root} of this rooted tree.
A {\it rooted tree-decomposition} of a graph $G$ is a tree-decomposition $(T,\X)$ of $G$ such that $T$ is a rooted tree.

Let $\F$ and $\F'$ be classes of graphs.
Let $\eta,\theta$ be nonnegative integers with $\eta \leq \theta$.
Recall that $\F^+=\F^{+1}$.
We say that a graph $G$ is {\it $(\eta,\theta,\F,\F')$-constructible} if there exists a rooted tree-decomposition $(T,\X)$ of $G$ of adhesion at most $\theta$ such that
	\begin{itemize}
		\item for every edge $tt' \in E(T)$, if $\lvert X_t \cap X_{t'} \rvert > \eta$, then one end of $tt'$ has no child, say $t'$, and the set $X_{t'}-X_t$ contains at most 1 vertex,
		\item for every $t \in V(T)$,
			\begin{itemize}
				\item if $t$ is the root of $T$, then $\lvert X_t \rvert \leq \theta$, 
				\item if $\eta>0$ and $t$ is the root of $T$, then $X_t \neq \emptyset$,
				\item if $t$ has a child in $T$, then $G[X_t] \in \F$,
				\item if $t$ has no child in $T$, then $G[X_t] \in \F^+$, and
				\item $\F'$ contains every graph that can be obtained from $G[X_t]$ by for each child $t'$ of $t$, adding a set $S_t$ of new vertices and adding new edges incident with $S_t$ such that the neighbors of each vertex in $S_t$ are contained in $X_t \cap X_{t'}$.
			\end{itemize}
	\end{itemize}
In this case, we call $(T,\X)$ an {\it $(\eta,\theta,\F,\F')$-construction} of $G$.

For every rooted tree $T$, define $I(T)$ to be the set of nodes of $T$ that have children.

Let $G$ be a graph and $m$ a positive integer.
Let $S \subseteq V(G)$.
Let $c: S \rightarrow [m]$ be a function.
Let $c'$ be an $m$-coloring of $G$ such that $c'(v)=c(v)$ for every $v \in S$.
Then we say that $c$ can be {\it extended} to $c'$.

Recall that a class $\F$ of graphs is hereditary if for every $G \in \F$, every induced subgraph of $G$ belongs to $\F$.
Note that if $\F$ is hereditary, then so is $\F^+$.

\begin{lemma} \label{tree_extension}
For any positive integers $\ell,N,m$ with $m \geq 2$, every nonnegative integer $\theta$ and any $(m,\ell,N)$-nice hereditary classes $\F,\F'$, there exists a function $f^*: {\mathbb N} \cup \{0\} \rightarrow {\mathbb N}$ such that the following holds. 
Let $\eta$ be a nonnegative integer with $\eta \leq \theta$.
Let $G$ be an $(\eta,\theta,\F,\F')$-constructible graph with an $(\eta,\theta,\F,\F')$-construction $(T,\X)$. 
Denote $\X$ by $(X_t: t \in V(T))$.
Let $t^*$ be the root of $T$.
Let $Z \subseteq N_G^{\leq 3\ell}[X_{t^*}]$.
If $c: Z \rightarrow [m]$ is a function, then $c$ can be extended to an $m$-coloring of $G^\ell$ with weak diameter in $G^\ell$ at most $f^*(\eta)$.
\end{lemma}

\begin{pf}
Let $\ell,N,m$ be positive integers with $m \geq 2$, and let $\theta$ be a nonnegative integer.
By Lemma \ref{apex_extension}, there exists an integer $N_{\F^+}$ such that $\F^+$ is $(m,\ell,N_{\F^+})$-nice.
Note that $\F \subseteq \F^+$, so we may assume that $N_{\F^+} \geq N$.
We define the following.
	\begin{itemize}
		\item Let $f_1: {\mathbb N} \rightarrow {\mathbb N}$ be the function such that for every $x \in {\mathbb N}$, $f_1(x)$ is the integer $N^*$ mentioned in Lemma \ref{deleting_centered_set} by taking $(k,r,\ell,N)=(\theta,3\ell,\ell,x)$. 
		\item Let $N_\theta$ be the integer $N^*$ mentioned in Lemma \ref{deleting_centered_set} by taking $(k,r,\ell,N)=(\theta,0,\ell,1)$.
			Note that we may assume that $N_\theta \geq \theta+1$.
		\item Let $N'_\theta$ be the number $N$ mentioned in Lemma \ref{all_centered} by taking $(k,r,\ell)=(\theta,3\ell,\ell)$. 
		\item Define $f^*: ({\mathbb N} \cup \{0\}) \rightarrow {\mathbb N}$ to be the function such that $f^*(0)=N_{\F^+} +N'_\theta+N_\theta+f_1(N)$, and for every $x \in {\mathbb N}$, $f^*(x)= \max\{(14\theta+4)\ell + 7\theta\ell^2f_1(f^*(x-1)),f^*(0)\}$.  
	\end{itemize}
Note that $f^*$ is an increasing function by Lemma \ref{deleting_centered_set}.

Let $\eta,G,(T,\X),t^*,Z,c$ be the ones as defined in the lemma.
Suppose to the contrary that $c$ cannot be extended to an $m$-coloring of $G^\ell$ with weak diameter in $G^\ell$ at most $f^*(\eta)$, and subject to this, the tuple $(\eta,\lvert I(T) \rvert+\lvert V(G)-Z \rvert+\lvert V(G) \rvert)$ is lexicographically minimal.

Since $(T,\X)$ is an $(\eta,\theta,\F,\F')$-construction, $\lvert X_{t^*} \rvert \leq \theta$.
So $Z$ is $(\theta,3\ell)$-centered.
If $Z=V(G)$, then $c$ can be extended to an $m$-coloring of $G^\ell$ with weak diameter in $G^\ell$ at most $N'_\theta \leq f^*(\eta)$ by Lemma \ref{all_centered}, a contradiction. 

So $Z \neq V(G)$.

\noindent{\bf Claim 1:} $\eta \geq 1$.

\noindent{\bf Proof of Claim 1:}
Suppose to the contrary that $\eta=0$.
Let $W=\{tt' \in E(T): X_t \cap X_{t'}=\emptyset\}$.
Note that for every component $C$ of $G$, $V(C) \subseteq X_{C'}$ for some component $C'$ of $T-W$.
Since every component of $G^\ell$ is $C^\ell$ for some component $C$ of $G$, there exists a component $C$ of $G$ such that $c|_{Z \cap V(C)}$ cannot be extended to an $m$-coloring of $C^\ell$ with weak diameter in $C^\ell$ at most $f^*(\eta)$.
Let $T_C$ be the component of $T-W$ with $V(C) \subseteq X_{T_C}$.
Let $t_C$ be the root of $T_C$.
Since $(T,\X)$ is a $(0,\theta,\F,\F')$-construction, $T_C$ is a star, and $G[X_{T_C}]$ is obtained from $G[X_{t_C}]$ by for each child $t'$ of $t_C$ in $T_C$, adding a vertex adjacent to some vertices in $X_{t_C} \cap X_{t'}$.
So $G[X_{T_C}] \in \F'$.
Let $G_C = G[X_{T_C}]-(Z \cap X_{T_C})$.
Since $\F'$ is hereditary, $G_C \in \F'$, so $G^\ell_C$ is $m$-colorable with weak diameter in $G_C^\ell$ at most $N \leq N_{\F^+} \leq f^*(\eta)$.
Since $\eta=0$ and $Z \subseteq N_G^{\leq 3\ell}[X_{t^*}]$, if $t_C \neq t^*$, then $Z \cap V(C) \subseteq Z \cap X_{T_C}=\emptyset$, so $c|_{Z \cap V(C)}$ can be extended to an $m$-coloring of $C^\ell$ with weak diameter in $C^\ell$ at most $f^*(\eta)$, a contradiction.
Hence $t_C=t^*$.
So $Z \cap X_{T_C} = Z$ is $(\theta,3\ell)$-centered in $C$.
By Lemma \ref{deleting_centered_set}, $c|_{Z \cap V(C)}$ can be extended to an $m$-coloring of $G_C^\ell$ with weak diameter in $G_C^\ell$ at most $f_1(N) \leq f^*(\eta)$, a contradiction.
$\Box$

\noindent{\bf Claim 2:} $G$ is connected.

\noindent{\bf Proof of Claim 2:}
Suppose to the contrary that $G$ is not connected.
Since every component of $G^\ell$ is $C^\ell$ for some component $C$ of $G$, there exists a component $C$ of $G$ such that $c|_{Z \cap V(C)}$ cannot be extended to an $m$-coloring of $C^\ell$ with weak diameter in $C^\ell$ at most $f^*(\eta)$.
Let $T_C$ be the subtree of $T$ induced by $\{t \in V(T): X_t \cap V(C) \neq \emptyset\}$.
Let $\X_C = (X_t \cap V(C): t \in V(T_C))$.

If $t^* \in V(T_C)$, then $t^*$ is the root of $T_C$ and $X_{t^*} \cap V(C) \neq \emptyset$, so $(T_C,\X_C)$ is an $(\eta,\theta,\F,\F')$-construction with $(\eta,\lvert I(T_C) \rvert+\lvert V(C)-(Z \cap V(C)) \rvert+\lvert V(C) \rvert)$ lexicographically smaller than $(\eta,\lvert I(T) \rvert+\lvert V(G)-Z \rvert+\lvert V(G) \rvert)$, a contradiction.
So $t^* \not \in V(T_C)$.
Since $Z \subseteq N_G^{\leq 3\ell}[X_{t^*}]$, $Z \cap V(C)=\emptyset$.
Let $T_C'$ be the rooted tree obtained from $T_C$ by adding a node $t^*_C$ adjacent to the root of $T_C$, where $t^*_C$ is the root of $T_C'$.
Let the bag at $t^*_C$ be the set consisting of a single vertex in the intersection of $V(C)$ and the bag of the root of $T_C$.
Then since $\eta \geq 1$ by Claim 1, we obtain an $(\eta,\theta,\F,\F')$-construction of $C$ with underlying tree $T_C'$.
Since $t^* \not \in V(T_C)$, $(\eta,\lvert I(T_C') \rvert+\lvert V(C)-(Z \cap V(C)) \rvert+\lvert V(C) \rvert)$ is lexicographically smaller than $(\eta,\lvert I(T) \rvert+\lvert V(G)-Z \rvert+\lvert V(G) \rvert)$.
Since $Z \cap V(C)=\emptyset$, the minimality implies that $c|_{Z \cap V(C)}$ can be extended to an $m$-coloring of $C^\ell$ with weak diameter in $C^\ell$ at most $f^*(\eta)$, a contradiction.
$\Box$

\noindent{\bf Claim 3:} $Z=N_G^{\leq 3\ell}[X_{t^*}]$ and $Z-X_{t^*} \neq \emptyset$.

\noindent{\bf Proof of Claim 3:}
Suppose to the contrary that there exists $v \in N_G^{\leq 3\ell}[X_{t^*}]-Z$.
Let $Z' = Z \cup \{v\}$.
Let $c': Z' \rightarrow [m]$ be the function obtained from $c$ by further defining $c'(v)=m$.
Then the minimality of $(\eta,\lvert I(T) \rvert+\lvert V(G)-Z \rvert+\lvert V(G) \rvert)$ implies that $c'$ (and hence $c$) can be extended to an $m$-coloring of $G^\ell$ with weak diameter in $G^\ell$ at most $f^*(\eta)$, a contradiction.

So $Z=N_G^{\leq 3\ell}[X_{t^*}]$.
Suppose to the contrary that $Z \subseteq X_{t^*}$.
So $N_G^{\leq 3\ell}[X_{t^*}]= Z \subseteq X_{t^*}$.
Since $G$ is connected by Claim 2, $V(G)=X_{t^*}= Z$, a contradiction.
$\Box$

Note that Claim 3 implies that $Z \neq \emptyset$.
And by Claim 1 and the definition of $(\eta,\theta,\F,\F')$-constructions, $X_{t^*} \neq \emptyset$.

For each $v \in X_{t^*}$, let $T_v$ be the subgraph of $T$ induced by $\{t \in V(T): N_G^{\leq 3\ell}[\{v\}] \cap X_t \neq \emptyset\}$.
Since $(T,\X)$ is a tree-decomposition, $T_v$ is a subtree of $T$ containing $t^*$.
So $\bigcup_{v \in X_{t^*}}T_v$ is a subtree of $T$ containing $t^*$.

Let $T_0 = \bigcup_{v \in X_{t^*}}T_v$.
Let $U_E = \{e \in E(T):$ exactly one end of $e$ is in $V(T_0)\}$. 
Let $U=\{t \in V(T_0): t$ is an end of $e$ for some $e \in U_E\}$. 
So for every $t \in U$, $t$ has a child, $t \in V(T_0)$, and for every $e_t \in U_E$ incident with $t$, the component of $T-e_t$ disjoint from $t$ is disjoint from $V(T_0)$.

For each $e \in E(T)$, define $X_e$ to be the intersection of the bags of the ends of $e$; define $T_e$ to be the component of $T-e$ disjoint from $t^*$.
Since the adhesion of $(T,\X)$ is at most $\theta$, $\lvert X_e \rvert \leq \theta$ for every $e \in E(T)$.

For each $e \in U_E$, define $\preceq_e$ to be the binary relation on $X_e$ such that for any $x,y \in X_e$, $x \preceq_e y$ if and only if there exists a path in $G[X_{T_e}]$ from $x$ to $y$ of length at most $7\ell$.
Clearly, for each $e \in U_E$, $\preceq_e$ is reflexive and symmetric, and since $\lvert X_e \rvert \leq \theta$, there exists a partition $\P_e$ of $X_e$ such that two vertices $x,y$ in $X_e$ are contained in the same part in $\P_e$ if and only if there exists a sequence $a_1,...,a_\theta$ (with not necessarily distinct terms) such that $a_1=x$, $a_\theta=y$, and $a_{i} \preceq_e a_{i+1}$ for every $i \in [\theta-1]$.

Let $Z_0 = N_G^{\leq 3\ell}[X_{t^*}]$.
So $Z_0=Z$ by Claim 3.
And for every $t \in V(T_0)$, $X_t \cap Z_0 \neq \emptyset$; for every $t \in V(T)-V(T_0)$, $X_t \cap Z_0=\emptyset$.

Define $G_0$ to be the graph obtained from $G[X_{T_0}]$ by for each $e \in U_E$ and $Y \in \P_e$, adding a new vertex $v_Y$ and edges incident with $v_Y$ such that $v_Y$ is adjacent to every vertex in $Y$. 
Note that $Z =Z_0 \subseteq V(G[X_{T_0}]) \subseteq V(G_0)$.

\noindent{\bf Claim 4:} There exists an $m$-coloring $c_0$ of $G_0^\ell$ with weak diameter in $G_0^\ell$ at most $f_1(f^*(\eta-1))$ such that $c_0(v)=c(v)$ for every $v \in Z$. 

\noindent{\bf Proof of Claim 4:}
Define $T'$ to be the rooted tree obtained from $T_0$ by for each $e \in U_E$ and $Y \in \P_e$, adding a node $t_Y$ adjacent to the end of $e$ in $V(T_0)$.
For each $t' \in V(T_0)$, define $X'_{t'}=X_{t'}$; for each $t' \in V(T')-V(T_0)$, $t'=t_Y$ for some $e \in U_E$ and $Y \in \P_e$, and we define $X'_{t'}=X_e \cup \{v_Y\}$. 
Let $\X'=(X'_t: t \in V(T'))$.

Clearly, $(T',\X')$ is a tree-decomposition of $G_0$ of adhesion at most $\max_{e \in U_E}\{\theta,\lvert X_e \rvert\} = \theta$.
For each $tt' \in E(T')$, say $t'$ is a child of $t$, if $tt' \in E(T_0)$, then $X'_t=X_t$, $X'_{t'}=X_{t'}$, $t$ has a child in both $T$ and $T'$, and $t'$ has a child in $T'$ if and only if $t'$ has a child in $T$; if $tt' \not \in E(T_0)$, then $t \in V(T_0)$ and $t' \not \in V(T_0)$, and $\lvert X'_{t'}-X'_t \rvert = 1$.
Hence for every $tt' \in E(T')$, if $\lvert X'_t \cap X'_{t'} \rvert > \eta$, then one end of $tt'$, say $t'$, has no child, and $\lvert X'_{t'}-X'_t \rvert \leq 1$.

Furthermore, $t^* \in V(T_0) \subseteq V(T')$ and $X'_{t^*}=X_{t^*}$, so $\lvert X'_{t^*} \rvert = \lvert X_{t^*} \rvert \leq \theta$.
Since $\eta \geq 1$ by Claim 1, $X'_{t^*} = X_{t^*} \neq \emptyset$.
In addition, for every $t \in V(T')$, if $t$ has a child in $T'$, then $t \in V(T_0) \subseteq V(T)$ has a child in $T$, so $G[X'_t]=G[X_t] \in \F$; if $t$ has no child in $T'$, then either $t \in V(T)$ has no child in $T$ (so $G[X'_t] = G[X_t] \in \F^+$), or $t \in V(T')-V(T)$ and $G[X'_t]$ can be obtained by adding a vertex to $G[X_e] \in \F$ for some $e \in U_E$, so $G[X'_t] \in \F^+$.
And if $t$ has a child in $T'$, then $t \in V(T_0) \subseteq V(T)$, so every graph that can be obtained from $G[X'_t]=G[X_t]$ by for each child $t'$ of $t$ in $T'$, adding a set $S_t$ of new vertices and adding new edges incident with $S_t$ such that the neighbors of each vertex in $S_t$ are contained in $X'_t \cap X'_{t'}$ belongs to $\F'$.

Therefore, $(T',\X')$ is an $(\eta,\theta,\F,\F')$-construction of $G_0$.

For every $t \in V(T')$, let $X''_t = X'_t-Z_0$.
Let $\X''=(X''_t: t\in V(T'))$.
So $(T',\X'')$ is a tree-decomposition of $G_0-Z_0$ of adhesion at most $\theta$.
Since for every $t \in V(T_0)$, $X_t \cap Z_0 \neq \emptyset$, we know for every $e \in E(T_0)$, $X_e \cap Z_0 \neq \emptyset$.
Note that $X''_{t^*} = X_{t^*}-Z_0=\emptyset$.
If $\eta-1=0$, then let $T'''=T'$ and $\X'''=\X''$; otherwise, $\eta-1 \geq 1$ by Claim 1, let $t_0$ be a node of $T'$ with $X''_{t_0} \neq \emptyset$ closest to $t^*$, let $v_0$ be a vertex in $X''_{t_0}$, let $T'''$ be the rooted tree obtained from $T'$ by adding a new node $t_0^*$ adjacent to $t^*$, where $t_0^*$ is the root of $T'''$, and let $\X'''=(X'''_t: t \in V(T'''))$, where $X'''_{t_0^*}=\{v_0\}$, $X'''_t=X''_t \cup \{v_0\}$ if $t \neq t_0^*$ and $t$ is in the path in $T'$ between $t^*$ and $t_0$, and $X'''_{t}=X''_t$ otherwise.

Then $(T''',\X''')$ is a tree-decomposition of $G_0-Z_0$.
Since $(T',\X')$ is an $(\eta,\theta,\F,\F')$-construction of $G_0$, and $G$ is connected, and $\F,\F'$ and $\F^+$ are hereditary, $(T''',\X''')$ is an $(\eta-1,\theta,\F,\F')$-construction of $G_0-Z_0$.

By the minimality of $\eta$, there exists an $m$-coloring $c_0'$ of $(G_0-Z_0)^\ell$ with weak diameter in $(G_0-Z_0)^\ell$ at most $f^*(\eta-1)$.
Let $c_0=c \cup c_0'$.
By Lemma \ref{deleting_centered_set}, $c_0$ is an $m$-coloring of $G_0^\ell$ with weak diameter in $G_0^\ell$ at most $f_1(f^*(\eta-1))$ with $c_0(v)=c(v)$ for every $v \in Z$.
$\Box$

Let $G_1 = G[\bigcup_{e \in U_E}X_{T_e}]$.
By the definition of $T_0$ and $U_E$, $Z \cap V(G_1)=\emptyset$. 
Let $S_1 = \bigcup_{e \in U_E}X_e$.
For $i \in [3]$, let $Z_i = N_{G_1}^{\leq i\ell}[S_1]$. 

\noindent{\bf Claim 5:} For every $e \in U_E$, $i \in [3]$ and $v \in Z_i \cap X_{T_e}-X_e$, there exists $Y \in \P_e$ such that $v \in N_{G[X_{T_e}]}^{\leq i\ell}[Y]$, and for every $Y' \in \P_e-\{Y\}$, $v \not \in N_{G[X_{T_e}]}^{\leq i\ell}[Y']$.

\noindent{\bf Proof of Claim 5:}
Since $v \in Z_{i} \cap X_{T_e}-X_e$, there exists a path $P$ in $G_1$ from $v$ to $X_e$ internally disjoint from $X_e$ of length at most $i\ell$.
Since $P$ is internally disjoint from $X_e$, $P$ is a path in $G[X_{T_e}]$.
Let $y$ be the vertex in $V(P) \cap X_e$.
Let $Y$ be the member of $\P_e$ containing $y$.
So $v \in N_{G[X_{T_e}]}^{\leq i\ell}[Y]$.
Let $Y'$ be any member of $\P_e-\{Y\}$.
If $v \in N_{G[X_{T_e}]}^{\leq i\ell}[Y']$, then there exists a walk in $G[X_{T_e}]$ from $Y$ to $Y'$ of length at most $2i\ell \leq 6\ell$, so $Y=Y'$ by the definition of $\P_e$, a contradiction.
Hence $v \not \in N_{G[X_{T_e}]}^{\leq i\ell}[Y']$.
$\Box$

For every $v \in Z_1-V(G_0) = Z_1-S_1 \subseteq Z_{3}-S_1$, there uniquely exists $e \in U_E$ such that $v \in X_{T_{e}}$, so by Claim 5, there exists a unique pair $(e_v,Y_v)$ with $e_v \in U_E$ and $Y_v \in \P_{e_v}$ such that $v \in N_{G[X_{T_{e_v}}]}^{\leq 3\ell}[Y_v]$.

Let $c_1: Z_1 \rightarrow [m]$ be the function such that 
	\begin{itemize}
		\item $c_1(u)=c_0(u)$ for every $u \in Z_1 \cap V(G_0) = S_1$, and 
		\item $c_1(u) = c_0(v_{Y_u})$ for every $u \in Z_1-V(G_0) = Z_1-S_1$.
	\end{itemize}
For each $i \in \{2,3\}$, let $c_i: Z_i \rightarrow [m]$ be the function such that 
	\begin{itemize}
		\item $c_i(v)=c_{i-1}(v)$ for every $v \in Z_{i-1}$, and 
		\item $c_i(v)=i-1$ for every $v \in Z_i-Z_{i-1}$.
	\end{itemize}

For every $e \in U_E$, let $G_e = G[X_{T_e}]$ and $Z_e = Z_3 \cap V(G_e)$, and let $c_e: Z_e \rightarrow [m]$ such that $c_e(v)=c_3(v)$ for every $v \in Z_e$.
By Claim 5, for every $e \in U_E$, $Z_e \subseteq N_{G_e}^{\leq 3\ell}[X_e]$.

\noindent{\bf Claim 6:} For every $e \in U_E$, $c_e$ can be extended to an $m$-coloring $c_e'$ of $G_e^\ell$ with weak diameter in $G_e^\ell$ at most $f^*(\eta)$.

\noindent{\bf Proof of Claim 6:}
Let $e \in U_E$.
If $\lvert X_e \rvert > \eta$, then since $(T,\X)$ is an $(\eta, \theta,\F,\F')$-construction of $G$, $\lvert V(G_e) \rvert \leq \lvert X_e \rvert+1 \leq \theta+1$, so $c_e$ can be extended to an $m$-coloring of $G_e^\ell$ with weak diameter in $G_e^\ell$ at most $\lvert V(G_e) \rvert \leq \theta+1 \leq f^*(\eta)$.

So we may assume that $\lvert X_e \rvert \leq \eta$.
If $X_e = \emptyset$, then $V(G_e)=\emptyset$ since $G$ is connected and $X_{t^*} \neq \emptyset$ by Claims 1 and 2, so we are done.

So we may assume that $X_e \neq \emptyset$.
Define $Q_e$ to be the rooted tree obtained from $T_e$ by adding a node $r_e$ adjacent to the end of $e$ in $V(T_e)$, where $r_e$ is the root of $Q_e$.
Let $W_{r_e} = X_e$; for every $t \in V(T_e)$, let $W_t=X_t$.
Let $\W=(W_t: t \in V(Q_e))$.
Then $(Q_e,\W)$ is a rooted tree-decomposition of $G_e$ of adhesion at most $\theta$ such that $\lvert W_{r_e} \rvert = \lvert X_e \rvert \leq \eta$.
So if $tt' \in E(Q_e)$ with $\lvert W_t \cap W_{t'} \rvert > \eta$, then $tt' \in E(T_e)$, so $W_t=X_t$ and $W_{t'}=X_{t'}$.
Since $G[W_{r_e}]=G[X_e]$ and $\F$ and $\F'$ are hereditary, $(Q_e,\W)$ is an $(\eta,\theta,\F,\F')$-construction of $G_e$.

Note that $I(Q_e) = \{r_e\} \cup (I(T) \cap V(T_e))$.
Note that every vertex that belongs to the shortest directed path in $T$ containing $t^*$ and an end of $e$ belongs to $I(T)-V(T_e)$.
So $\lvert I(Q_e) \rvert \leq \lvert I(T) \rvert$, and equality holds only when $t^*$ is an end of $e$. 
If $t^*$ is an end of $e$, then since $X_e \neq \emptyset$, by Claim 3, $X_t \cap Z \supseteq X_t \cap X_{t^*} \neq \emptyset$, where $t$ is the end of $e$ other than $t^*$, so $t \in V(T_0)$, a contradiction.
Hence $\lvert I(Q_e) \rvert < \lvert I(T) \rvert$.

Recall that $Z_e \subseteq N_{G_e}^{\leq 3\ell}[X_e] = N_{G_e}^{\leq 3\ell}[W_{r_e}]$.
So by the minimality of $(\eta,\lvert I(T) \rvert+\lvert V(G)-Z \rvert+\lvert V(G) \rvert)$, $c_e$ can be extended to an $m$-coloring of $G_e^\ell$ with weak diameter in $G_e^\ell$ at most $f^*(\eta)$.
$\Box$

Let $c_0$ be the $m$-coloring mentioned in Claim 4.
For every $e \in U_E$, let $c_e'$ be the $m$-coloring mentioned in Claim 6.

Define $c^* = c_0|_{V(G) \cap V(G_0)} \cup \bigcup_{e \in U_E}c'_e$.
Note that $c^*$ is well-defined by the definition of $c_1$, and $c^*$ is an $m$-coloring of $G^\ell$ that can be obtained from extending $c$ by Claim 4.
So there exists a $c^*$-monochromatic component $M$ in $G^\ell$ with weak diameter in $G^\ell$ greater than $f^*(\eta)$.

\noindent{\bf Claim 7:} There exists no $e \in U_E$ such that $V(M) \cap (Z_1-X_e) \cap V(G_e) \neq \emptyset$ and $V(M) \cap (Z_3-Z_2) \cap V(G_e) \neq \emptyset$.

\noindent{\bf Proof of Claim 7:}
Suppose to the contrary that there exists $e \in U_E$ such that $V(M) \cap (Z_1-X_e) \cap V(G_e) \neq \emptyset$ and $V(M) \cap (Z_3-Z_2) \cap V(G_e) \neq \emptyset$.
Then there exists a path $P$ in $M \subseteq G^\ell$ from $V(M) \cap (Z_1-X_e) \cap V(G_e)$ to $V(M) \cap (Z_3-Z_2) \cap V(G_e)$.
Hence $P$ contains an edge $e_1$ between $V(M) \cap (Z_1-X_e) \cap V(G_e)$ and $V(M) \cap (Z_2-Z_1) \cap V(G_e)$, and $P$ contains an edge $e_2$ between $V(M) \cap (Z_2-Z_1) \cap V(G_e)$ and $V(M) \cap (Z_3-Z_2) \cap V(G_e)$.
In particular, $M$ contains a vertex $v_1$ in $(Z_2-Z_1) \cap V(G_e)$ and a vertex $v_2$ in $(Z_3-Z_2) \cap V(G_e)$.
By the definition of $c_e,c_2,c_3$, $c^*(v_1)=c_2(v_1)=1$ and $c^*(v_2)=c_3(v_2)=2$.
So $M$ is not $c^*$-monochromatic, a contradiction.
$\Box$

\noindent{\bf Claim 8:} For every $e \in U_E$, if $V(M) \cap (Z_1-X_e) \cap V(G_e)= \emptyset$, then $V(M) \cap V(G_e)-X_e = \emptyset$.

\noindent{\bf Proof of Claim 8:}
Suppose to the contrary that there exists $e \in U_E$ such that $V(M) \cap (Z_1-X_e) \cap V(G_e)= \emptyset$ and $V(M) \cap V(G_e)-X_e \neq \emptyset$.
Then $V(M) \subseteq V(G_e)-X_e$ since $M \subseteq G^\ell$.
Since $V(M) \cap (Z_1-X_e) \cap V(G_e)= \emptyset$, $V(M) \subseteq V(G_e)-Z_1$.
Hence every edge in $M \subseteq G^\ell$ is an edge in $G_e^\ell$.
So $M$ is a $c'_e$-monochromatic component in $G_e^\ell$.
Hence the weak diameter of $M$ in $G_e^\ell$ is at most $f^*(\eta)$ by Claim 6.
That is, for any vertices $x,y$ in $M$, there exists a path $P_{x,y}$ in $G_e^\ell$ between $x$ and $y$ of length at most $f^*(\eta)$.
Since $G_e \subseteq G$, $G_e^\ell \subseteq G^\ell$.
So for any vertices $x,y$ in $M$, $P_{x,y}$ is a path in $G^\ell$ between $x$ and $y$ of length at most $f^*(\eta)$.
Hence the weak diameter of $M$ in $G^\ell$ is at most $f^*(\eta)$, a contradiction.
$\Box$

Let $B = \{e \in U_E: V(M) \cap V(G_e)-X_e \neq \emptyset\}$.

\noindent{\bf Claim 9:} $V(M) \subseteq (V(G_0) \cap V(G)) \cup \bigcup_{e \in B}(X_{T_e} \cap Z_2)$.

\noindent{\bf Proof of Claim 9:}
For every $e \in B$, by Claim 8, $V(M) \cap (Z_1-X_e) \cap V(G_e) \neq \emptyset$, so by Claim 7, $V(M) \cap (Z_3-Z_2) \cap V(G_e)=\emptyset$.
Hence for every $e \in B$, since $M \subseteq G^\ell$, $V(M) \cap V(G_e) \subseteq V(G_e) \cap Z_2 = X_{T_e} \cap Z_2$.
Therefore, $V(M) \subseteq (V(G_0) \cap V(G)) \cup \bigcup_{e \in B}(X_{T_e} \cap Z_2)$.
$\Box$

\noindent{\bf Claim 10:} For every $x \in V(M)-(V(G_0) \cap V(G))$, there exists $e \in B$ with $x \in X_{T_e} \cap Z_2-X_e$, there uniquely exists $Y_x \in \P_e$ such that $x \in N_{G[X_{T_e}]}^{\leq 2\ell}[Y_x]$, and $c^*(M)=c_0(v_{Y_x})$.

\noindent{\bf Proof of Claim 10:}
Let $x \in V(M)-(V(G_0) \cap V(G))$.
By Claim 9, there exists $e \in B$ such that $x \in X_{T_e} \cap Z_2-X_e$.
So by Claim 5, there uniquely exists $Y_x \in \P_e$ such that $x \in N_{G[X_{T_e}]}^{\leq 2\ell}[Y_x]$.
Furthermore, $V(M) \cap (Z_1-X_e) \cap V(G_e) \neq \emptyset$ by Claim 8.
So there exists $x' \in V(M) \cap (Z_1-X_e) \cap V(G_e)$ such that $c^*(x)=c^*(x')=c^*(M)$.

We further choose such $x'$ such that the distance in $M$ between $x$ and $x'$ is as small as possible.
(Note that it is possible that $x=x'$.)
Let $P$ be a shortest path in $M \subseteq G^\ell$ between $x$ and $x'$.
By the choice of $x'$, $V(P)-\{x'\} \subseteq V(M) \cap X_{T_e} \cap Z_2-Z_1$.
Denote $P$ by $x_1x_2...x_{\lvert V(P) \rvert}$. 
So for every $i \in [\lvert V(P) \rvert-1]$, there exists a path in $G_e$ between $x_i$ and $x_{i+1}$ with length in $G_e$ at most $\ell$. 

By Claim 5, there uniquely exists $Y_{x'} \in \P_e$ such that $x' \in N_{G[X_{T_e}]}^{\leq \ell}[Y_{x'}]$, and for each $v \in V(P)-\{x'\}$, there uniquely exists $Y_{v} \in \P_e$ such that $v \in N_{G[X_{T_e}]}^{\leq 2\ell}[Y_{v}]$.
For any $i \in [\lvert V(P) \rvert-1]$, if $Y_{x_i} \neq Y_{x_{i+1}}$, then there exists a walk in $G_e$ from $Y_{x_i}$ to $Y_{x_{i+1}}$ with length at most $2\ell+\ell+2\ell \leq 5\ell$, so $Y_{x_i}=Y_{x_{i+1}}$ by the definition of $\P_e$, a contradiction.
So $Y_{x_i}=Y_{x_{i+1}}$ for every $i \in [\lvert V(P) \rvert-1]$.
In particular, $v_{Y_x}=v_{Y_{x'}}$.
Hence by the definition of $c_1$, $c^*(M)=c^*(x')=c_1(x')=c_0(v_{Y_{x'}})=c_0(v_{Y_{x}})$.
$\Box$

For every $x \in V(M)-(V(G_0) \cap V(G))$, define $Y_x$ to be the set mentioned in Claim 10.
Let $M'$ be the graph obtained from $M$ by for each $e \in B$ and $Y \in \P_e$, identifying all vertices $x \in V(M)-(V(G_0) \cap V(G))$ with $Y_x=Y$ into a vertex $v_Y$.
Note that there is a natural injection from $V(M')$ to $V(G_0)$ obtained by the identification mentioned in the definition of $M'$.
So we may assume $V(M') \subseteq V(G_0)$.

\noindent{\bf Claim 11:} $M'$ is contained in a $c_0$-monochromatic component in $G_0^\ell$. 

\noindent{\bf Proof of Claim 11:}
By Claim 10, for every $x \in V(M)-(V(G_0) \cap V(G))$, $c^*(M)=c_0(v_{Y_x})$.
So all vertices in $M'$ have the same color in $c_0$.
Hence to prove that $M'$ is contained in a $c_0$-monochromatic component in $G_0^\ell$, it suffices to prove that $M'$ is a subgraph of $G_0^\ell$.
Since $V(M') \subseteq V(G_0)=V(G_0^\ell)$, it suffices to prove that $E(M') \subseteq E(G_0^\ell)$.

Note that for any $e \in U_E$, distinct vertices $x,y \in X_e$ and path $P$ in $G[X_{T_e}]$ between $x$ and $y$ internally disjoint from $X_e$ of length at most $\ell$ having at least one internal vertex, the ends of $P$ are contained in the same part (say $Y$) of $\P_e$, so there exists a path $\overline{P}=xv_Yy$ in $G_0$ of length two between $x$ and $y$; since $P$ has at least one internal vertex, the length of $\overline{P}$ is at most the length of $P$.
Hence, for every path $P$ in $G$ of length at most $\ell$ between two distinct vertices in $V(G_0) \cap V(G)$, we can replace each maximal subpath $P'$ of $P$ of length at least two whose all internal vertices are not in $V(G_0)$ by $\overline{P'}$ to obtain a walk $\overline{P}$ in $G_0$ of length at most the length of $P$ having the same ends as $P$.

Hence if $xy$ is an edge of $M'$ with $x,y \in V(G_0) \cap V(G)$, then $x,y \in V(M)$, and since $M \subseteq G^\ell$, there exists a path $P_{xy}$ in $G$ of length at most $\ell$ between $x$ and $y$, so $\overline{P_{xy}}$ is a walk in $G_0$ of length at most $\ell$ between $x$ and $y$, so $xy \in E(G_0^\ell)$.

Now assume that there exists an edge $xy$ of $M'$ with $x \in V(G_0) \cap V(G)$ and $y \not \in V(G_0) \cap V(G)$.
Since $y \not \in V(G_0) \cap V(G)$, there exists $y_0 \in V(M)$ with $v_{Y_{y_0}}=y$ such that $xy_0 \in E(M)$, and there exists $e_y \in U_E$ such that $y_0 \in X_{T_{e_y}}-X_{e_y}$.
Since $M \subseteq G^\ell$, there exists a path $P_{xy}$ in $G$ of length at most $\ell$ between $x$ and $y_0$.
Let $y'$ be the vertex in $V(P_{xy}) \cap X_{e_y}$ such that the subpath of $P_{xy}$ between $y_0$ and $y'$ is contained in $G[X_{T_{e_y}}]$.
Then $y' \in Y_{y_0}$ by Claim 5. 
So $x\overline{P'_{xy}}y'y$ is a walk in $G_0^\ell$ of length at most $\ell$, where $P'_{xy}$ is the subpath of $P_{xy}$ between $x$ and $y'$.
So $xy \in E(G_0^\ell)$.

Hence every edge of $M'$ incident with a vertex of $V(G_0) \cap V(G)$ is an edge of $G_0^\ell$.

Now assume that there exist $e \in B$ and distinct $Y,Y' \in \P_e$ such that $v_Yv_{Y'} \in E(M')$.
So there exists $ab \in E(M)$ such that $a,b \in V(M) \cap X_{T_e}-X_e$, $Y_a=Y$ and $Y_b=Y'$.
By Claim 10, there exist a path $P_a$ in $G[X_{T_e}]$ from $a$ to $Y_a$ with length at most $2\ell$ and a path $P_b$ in $G[X_{T_e}]$ from $b$ to $Y_b$ of length at most $2\ell$.
Since $ab \in E(M) \subseteq E(G^\ell)$, there exists a path $P_{ab}$ in $G$ of length at most $\ell$ from $a$ to $b$.
If $V(P_{ab}) \subseteq X_{T_e}$, then $P_a \cup P_{ab} \cup P_b$ is a walk in $G[X_{T_e}]$ from $Y_a$ to $Y_b$ of length at most $2\ell+\ell+2\ell<7\ell$, contradicting that $Y_a=Y$ and $Y_b=Y'$ are distinct parts of $\P_e$.
So $V(P_{ab}) \not \subseteq X_{T_e}$.
In particular, there exist distinct $a',b' \in V(P_{ab}) \cap X_e$ such that the subpath $P_{a'}$ of $P_{ab}$ between $a$ and $a'$ and the subpath $P_{b'}$ of $P_{ab}$ between $b$ and $b'$ are contained in $G[X_{T_e}]$.
Since the length of $P_{a'}$ and $P_{b'}$ are at most $\ell$, $a' \in Y_a$ and $b' \in Y_b$ by Claim 5.
So $v_Ya'\overline{P'}b'v_{Y'}$ is path in $G_0$ of length at most the length of $P_{ab}$, where $P'$ is the subpath of $P_{ab}$ between $a'$ and $b'$.
Hence $v_Yv_{Y'} \in G_0^\ell$.

Finally, assume that there exist distinct $e_1,e_2 \in B$, $Y_1 \in \P_{e_1}$ and $Y_2 \in \P_{e_2}$ such that $v_{Y_1}v_{Y_2} \in E(M')$.
So for each $i \in [2]$, there exists $x_i \in V(M) \cap Z_2 \cap X_{T_{e_i}}-X_{e_i}$ by Claim 10 such that $Y_{x_i}=Y_i$, and there exists a path $P_{x_1x_2}$ in $G$ of length at most $\ell$ between $x_1$ and $x_2$.
For each $i \in [2]$, let $y_i$ be the vertex in $V(P_{x_1x_2}) \cap X_{e_i}$ such that the subpath of $P_{x_1x_2}$ between $x_i$ and $y_i$ is contained in $G[X_{T_{e_i}}]$.
Then $v_{Y_1}y_1\overline{P'_{x_1x_2}}y_2v_{Y_2}$ is a walk in $G_0$ of length at most $\ell$, where $P'_{x_1x_2}$ is the subpath of $P_{x_1x_2}$ between $y_1$ and $y_2$.
Therefore, $v_{Y_1}v_{Y_2} \in E(G_0^\ell)$.

This proves $E(M') \subseteq E(G_0^\ell)$, and hence the claim follows.
$\Box$

\noindent{\bf Claim 12:} For any two vertices $x,y \in V(G) \cap V(G_0)$, if $P$ is a path in $G_0$ between $x$ and $y$, then there exists a walk $\hat{P}$ in $G$ between $x$ and $y$ of length at most $7\theta\ell\lvert E(P) \rvert$.

\noindent{\bf Proof of Claim 12:}
Let $Q$ be a path in $G_0$ such that there exists $e \in U_E$ such that $Q$ is from $X_e$ to $X_e$ internally disjoint from $X_e$ and contains at least one internal vertex corresponding to  a vertex in $X_{T_e}-X_e$.
Then $Q$ has length two, and there exists $Y_Q \in \P_e$ containing both ends of $Q$.
So there exists a path $\hat{Q}$ in $G[X_{T_e}]$ between the ends of $Q$ of length at most $7\theta\ell$ by the definition of $\P_e$.

Let $x,y \in V(G) \cap V(G_0)$.
Let $P$ be a path in $G_0$ between $x$ and $y$.
Then by replacing each subpath $P'$ of $P$ in which there exists $e_{P'} \in U_E$ such that $P'$ is from $X_{e_{P'}}$ to $X_{e_{P'}}$ internally disjoint from $X_{e_{P'}}$ and contains at least one internal vertex corresponding to a vertex in  $X_{T_{e_{P'}}}-X_{e_{P'}}$ by $\hat{P'}$, we obtain a walk in $G$ between $x$ and $y$ of length at most $7\theta\ell \lvert E(P) \rvert$.
$\Box$

By Claim 11, $M'$ is contained in a $c_0$-monochromatic component in $G_0^\ell$.
By Claim 4, the weak diameter of $M'$ in $G_0^\ell$ is at most $f_1(f^*(\eta-1))$.

For every vertex $x \in V(M)$, if $x \in V(G_0) \cap V(G)$, then define $h(x)=x$; otherwise, define $h(x)=v_{Y_x}$.
Note that for every $x \in V(M)$, $h(x) \in V(M')$.

\noindent{\bf Claim 13:} For any vertices $x,y \in V(M)$, there exists a path in $G$ between $x$ and $y$ with length at most $f^*(\eta)$. 

\noindent{\bf Proof of Claim 13:}
Since $h(x),h(y) \in V(M')$ and the weak diameter of $M'$ in $G_0^\ell$ is at most $f_1(f^*(\eta-1))$, there exists a path $P_0$ in $G_0^\ell$ between $h(x)$ and $h(y)$ of length at most $f_1(f^*(\eta-1))$.
So there exists a path $P_0'$ in $G_0$ of length at most $\ell \cdot f_1(f^*(\eta-1))$ between $h(x)$ and $h(y)$.

For each $u \in \{x,y\}$, if $h(u) \in V(G_0) \cap V(G)$, then $h(u)=u$, and we let $h'(u)=u$ and $h''(u)=u$; otherwise, by Claim 10, there exists $e \in U_E$ such that $u \in (X_{T_e}-X_e) \cap N_{G[X_{T_e}]}^{\leq 2\ell}[Y_u]$, so there exists a path $P_u$ in $G[X_{T_e}]$ from $u$ to $Y_u$ with length at most $2\ell$, and we let $h'(u)$ be the end of $P_u$ in $Y_u$, and let $h''(u)$ be the neighbor of $h(u)$ in $P_0'$.
Note that for the latter case, $h''(u) \in Y_u$ by the definition of $G_0$ and $M'$.
So for each $u \in \{x,y\}$, there exists a path in $G$ from $u$ to $h'(u)$ of length at most $2\ell$, and there exists a path in $G$ from $h'(u)$ to $h''(u)$ of length at most $7\theta\ell$ by the definition of $\P_e$, so there exists a walk $W_u$ in $G$ from $u$ to $h''(u)$ of length at most $(7\theta+2)\ell$.

Note that $h''(x)$ and $h''(y)$ are contained in $P_0'$.
So the subpath $P'$ of $P_0'$ between $h''(x)$ and $h''(y)$ is a path in $G_0$ of length at most $\ell \cdot f_1(f^*(\eta-1))$.
Since $h''(x), h''(y) \in V(G_0) \cap V(G)$, by Claim 12, there exists a walk $\hat{P'}$ in $G$ between $h''(x)$ and $h''(y)$ of length at most $7\theta\ell\lvert E(P') \rvert \leq 7\theta\ell^2f_1(f^*(\eta-1))$.
Therefore, $W_x \cup \hat{P'} \cup W_y$ is a walk in $G$ from $x$ to $y$ of length at most $2(7\theta+2)\ell + 7\theta\ell^2f_1(f^*(\eta-1)) \leq f^*(\eta)$.
$\Box$

Hence the weak diameter of $M$ in $G$ is at most $f^*(\eta)$ by Claim 13.
Therefore, the weak diameter of $M$ in $G^\ell$ is at most $f^*(\eta)$, a contradiction.
This proves the lemma.
\end{pf}

\bigskip

Now we are ready to prove Theorem \ref{tree_extension_intro}.
The following is a restatement.

\begin{theorem} \label{tree_extension_ad}
Let $\F$ and $\F'$ be hereditary classes of graphs.
Let $\theta$ be a positive integer.
Let $\F^*$ be a class of graphs such that for every $G \in \F^*$, there exists a tree-decomposition $(T,\X)$ of $G$ of adhesion at most $\theta$, where $\X=(X_t: t \in V(T))$, such that for every $t \in V(T)$, 
	\begin{itemize}
		\item $G[X_t] \in \F$, and 
		\item $\F'$ contains every graph that can be obtained from $G[X_t]$ by adding new vertices and new edges such that for each new vertex $v$, there exists a neighbor $t_v$ of $t$ in $T$ such that the neighbors of $v$ are contained in $X_t \cap X_{t_v}$.
	\end{itemize}
Then $\ad(\F^*) \leq \max\{\ad(\F), \ad(\F'),1\}$. 
\end{theorem}

\begin{pf}
By Proposition \ref{wd_color_ad}, there exists a function $f: {\mathbb N} \rightarrow {\mathbb N}$ such that $\F \cup \F'$ is $(\ad(\F \cup \F')+1,\ell,f(\ell))$-nice for every $\ell \in {\mathbb N}$.
Define $g: {\mathbb N} \rightarrow {\mathbb N}$ to be the function such that for every $x \in {\mathbb N}$, $g(x)=f^*_x(\theta)$, where $f^*_x$ is the function $f^*$ mentioned in Lemma \ref{tree_extension} by taking $(\ell,N,m,\theta,\F,\F')=(x,f(x),\max\{\ad(\F \cup \F'),1\}+1,\theta,\F,\F')$.

Let $G \in \F^*$.
So there exists a tree-decomposition $(T,\X)$ of $G$ of adhesion at most $\theta$, where $\X=(X_t: t \in V(T))$, such that for every $t \in V(T)$, $G[X_t] \in \F$, and $\F'$ contains every graph that can be obtained from $G[X_t]$ by adding new vertices and new edges such that for each new vertex $v$, there exists a neighbor $t_v$ of $t$ in $T$ such that the neighbors of $v$ are contained in $X_t \cap X_{t_v}$.
Let $t_0$ be a node of $T$ with $X_{t_0} \neq \emptyset$, and let $v_0$ be a vertex in $X_{t_0}$.
Let $T'$ be the rooted tree obtained from $T$ by adding a new node $t_0'$ adjacent to $t_0$, where $t_0'$ is the root of $T'$.
Let $X'_{t_0'}=\{v_0\}$; for every $t \in V(T)$, let $X'_t=X_t$.
Let $\X'=(X'_t: t \in V(T'))$.
Then $(T',\X')$ is a $(\theta,\theta,\F,\F')$-construction of $G$.
For every $\ell \in {\mathbb N}$, applying Lemma \ref{tree_extension} by taking $(\ell,N,m,\theta,\F,\F',\eta,Z)=(\ell,f(\ell),\max\{\ad(\F \cup \F'),1\}+1,\theta,\F,\F',\theta,\emptyset)$, $G^\ell$ is $(\max\{\ad(\F \cup \F'),1\}+1)$-colorable with weak diameter in $G^\ell$ at most $g(\ell)$.

Hence $\F^*$ is $(\max\{\ad(\F \cup \F'),1\}+1,\ell,g(\ell))$-nice for every $\ell \in {\mathbb N}$.
By Proposition \ref{wd_color_ad}, $\ad(\F^*) \leq \max\{\ad(\F \cup \F'),1\} = \max\{\ad(\F),\ad(\F'),1\}$.
\end{pf}

\section{Applications to asymptotic dimension} \label{sec:app_ad}

Now we are ready to prove Theorem \ref{tw_ad_intro_2}.
The following is a restatement.

\begin{theorem} \label{tw_ad}
Let $w$ be a positive integer.
Let $\F$ be the class of graphs of tree-width at most $w$.
Then $\ad(\F)=1$.
\end{theorem}

\begin{pf}
Let $\F_1$ be the class of graphs on at most $w+1$ vertices.
Let $\F_2$ be the class of graphs that have a vertex-cover of size at most $w+1$.
Note that $\F_1,\F_2$ are hereditary classes and $\F_1 \subseteq \F_2$.
So $\max\{\ad(\F_1), \ad(\F_2)\}=\ad(\F_2) \leq 0$ by Lemma \ref{vc_ad}.

Note that for every graph $G$ of tree-width at most $w$, there exists a tree-decomposition $(T,\X)$ of $G$ of adhesion at most $w+1$, where $\X=(X_t: t \in V(T))$, such that for every $t \in V(T)$, $G[X_t] \in \F_1$, and $\F_2$ contains every graph that can be obtained from $G[X_t]$ by adding new vertices and new edges such that for each new vertex $v$, there exists a neighbor $t_v$ of $t$ in $T$ such that the neighbors of $v$ are contained in $X_t \cap X_{t_v}$.
Hence by Theorem \ref{tree_extension_ad}, $\ad(\F) \leq \max\{\ad(\F_1), \ad(\F_2),1\}=1$.

Since every path has tree-width at most 1, and there exists no integer $N$ such that every path is 1-colorable with weak diameter at most $N$, $\ad(\F) \geq 1$.
This shows $\ad(\F)=1$.
\end{pf}

\bigskip

Note that Theorem \ref{tw_ad_intro_1} follows from Theorem \ref{tw_ad_intro_2} by the Grid Minor Theorem \cite{rs_V}.
To prove Theorem \ref{layered_tw_ad_intro}, we need the following.

\begin{theorem}[{\cite[Theorem 5.2]{bbegps}}] \label{layer_trick}
Let $n$ be an integer such that for every nonnegative integer $k$, the class of graphs of tree-width at most $k$ has asymptotic dimension at most $n$.
Then for every positive integer $w$, the class of graphs of layered tree-width at most $w$ has asymptotic dimension at most $n+1$.
\end{theorem}

The following is a restatement of Theorem \ref{layered_tw_ad_intro}.

\begin{theorem} \label{layered_tw_ad}
Let $w$ be a positive integer.
Then the class of graphs of layered tree-width at most $w$ is at most 2.
\end{theorem}

\begin{pf}
It immediately follows from Theorems \ref{tw_ad} and \ref{layer_trick}.
\end{pf}

\bigskip

We need the following lemma in order to prove Theorem \ref{minor_ad_intro}.

\begin{lemma} \label{small_exten_minor}
Let $p$ be a positive integer.
For every positive integer $x$, let $\F_x$ be the class of graphs of layered tree-width at most $x$.
Let $\W$ be the class of graphs such that for every $G \in \W$, $G$ can be obtained from a graph $G' \in \F_p^{+p}$ by adding new vertices and new edges incident with new vertices such that for each new vertex, its neighbors are contained in a clique in $G'$.
Then $\W \subseteq \F_{p+1}^{+p}$.
\end{lemma}

\begin{pf}
Let $G \in \W$.
So there exists $H \in \F_p^{+p}$ such that $G$ can be obtained from $H$ by adding new vertices and new edges incident with new vertices such that for each new vertex, its neighbors are contained in a clique in $H$.
Since $H \in \F_p^{+p}$, there exists $Z \subseteq V(H)$ with $\lvert Z \rvert \leq p$ such that $H-Z \in \F_p$.
Hence $G-Z$ can be obtained from $H-Z$ by adding new vertices and new edges incident with new vertices such that for each vertex $v \in V(G)-(V(H) \cup Z)=V(G)-V(H)$, its neighbors are contained in a clique $C_v$ in $H-Z$.

Since $H-Z \in \F_p$, there exist a layering $(V_1,V_2,...)$ of $H-Z$ and a tree-decomposition $(T,\X=(X_t: t \in V(T))$ of $H-Z$ such that the intersection of any $V_i$ and $X_t$ has size at most $p$.
For every $v \in V(G)-V(H)$, since $C_v$ is a clique in $H-Z$, there exist $t_{C_v} \in V(T)$ with $C_v \subseteq X_{t_{C_v}}$ and a positive integer $i_{C_v}$ such that $C_v \subseteq V_{i_{C_v}} \cup V_{i_{C_v}+1}$.

For each positive integer $i$, define $V_i' = V_i \cup \{v \in V(G)-V(H): i_{C_v}=i\}$.
Then $(V_1',V_2',...)$ is a layering of $G-Z$.
Let $T'$ be the tree obtained from $T$ by for every $v \in V(G)-V(H)$, adding a new node $t_v$ adjacent to $t_{C_v}$.
For every $t \in V(T)$, define $X'_t = X_t$; for every $t \in V(T')-V(T)$, $t=t_v$ for some $v \in V(G)-V(H)$, and we define $X'_t = C_v \cup \{v\}$.
Let $\X'=(X'_t: t \in V(T'))$.
Then $(T',\X')$ is a tree-decomposition of $G-Z$.

Let $i$ be a positive integer, and let $t \in V(T')$.
If $t \in V(T)$, then $X_t' \subseteq V(H)-Z$, so $\lvert X'_t \cap V'_i \rvert = \lvert X_t \cap V_i \rvert \leq p$.
If $t \in V(T')-V(T)$, then there exists $v \in V(G)-V(H)$ such that $t=t_v$, so $X'_t \cap V_i' = (C_v \cap V_i) \cup (\{v\} \cap V_i') \subseteq (X_{t_{C_v}} \cap V_i) \cup \{v\}$, and hence $\lvert X'_t \cap V_i' \rvert \leq p+1$.

Therefore, the layered tree-width of $G-Z$ is at most $p+1$.
So $G-Z \in \F_{p+1}$.
Hence $G \in \F_{p+1}^{+p}$.
This shows $\W \subseteq \F_{p+1}^{+p}$.
\end{pf}

\bigskip

Let $G$ be a graph.
Let $(T,\X)$ be a tree-decomposition of $G$, where $\X=(X_t: t \in V(T))$.
For every $t \in V(T)$, the {\it torso} at $t$ is the graph obtained from $G[X_t]$ by adding edges such that for every neighbor $t'$ of $t$ in $T$, $X_t \cap X_{t'}$ is a clique.

The following is a restatement of Theorem \ref{minor_ad_intro}.

\begin{theorem} \label{minor_ad}
Let $H$ be a graph.
Let $\F$ be the class of graphs with no $H$-minor.
Then $\ad(\F) \leq 2$.
\end{theorem}

\begin{pf}
For every positive integer $x$, let $\F_x$ be the class of graphs of layered tree-width at most $x$.
By \cite[Theorem 1.3]{rs_XVI} and \cite[Theorem 20]{dmw}, there exists a positive integer $p$ such that for every graph $G \in \F$, there exists a tree-decomposition $(T,\X)$ of $G$ of adhesion\footnote{We remark that the bounded adhesion conclusion is not explicitly stated in \cite[Theorem 1.3]{rs_XVI}, but it can be easily derived from it: the description of the torsos in the theorem implies that the maximum size of a clique in torsos is bounded, so the tree-decomposition has bounded adhesion.} at most $p$ such that for every $t \in V(T)$, the torso at $t$ belongs to $\F_p^{+p}$. 

Let $\W$ be the class of graphs such that for every $G \in \W$, $G$ can be obtained from a graph $G' \in \F_p^{+p}$ by adding new vertices and new edges incident with new vertices such that for each new vertex, its neighbors are contained in a clique in $G'$.
By Lemma \ref{small_exten_minor}, $\W \subseteq \F_{p+1}^{+p}$.

Note that $\F_p^{+p}$ and $\F_{p+1}^{+p}$ are closed under taking subgraphs.
Hence for every $G \in \F$, there exists a tree-decomposition $(T,\X)$ of $G$ of adhesion at most $p$, where $\X=(X_t: t \in V(T))$, such that for every $t \in V(T)$, $G[X_t] \in \F_p^{+p}$, and $\F_{p+1}^{+p}$ contains every graph that can be obtained from $G[X_t]$ by adding new vertices and new edges such that for each new vertex $v$, there exists a neighbor $t_v$ of $t$ in $T$ such that the neighbors of $v$ are contained in $X_t \cap X_{t_v}$.

Therefore, by Theorem \ref{tree_extension_ad}, $\ad(\F) \leq \max\{\ad(\F_p^{+p}), \ad(\F_{p+1}^{+p}),1\} \leq 2$, where the last inequality follows from Theorems \ref{apex_extension_ad} and \ref{layered_tw_ad}.
\end{pf}

\bigskip

\bigskip

\noindent{\bf Acknowledgement:}
The author thanks Louis Esperet for bringing \cite{bbegps} which is the main motivation of this paper to his attention and for some discussions.
The author also thanks him for suggesting a shorter proof of Lemma \ref{deleting_centered_set}, which is included in this version of the paper.
The author thanks David Wood for pointing out consequences of \cite{demww,dmw_2} and thanks David Hume for pointing out the implicitly proved result in \cite{bst} as stated in the introduction.
The author thanks anonymous reviewers for their careful reading and suggestions.

\end{document}